%% file: sn-article.tex
\documentclass[sn-mathphys,Numbered]{sn-jnl}

\usepackage{graphicx}%
\usepackage{multirow}%
\usepackage{amsmath,amssymb,amsfonts,nicefrac}%
\usepackage{amsthm}%
\usepackage{mathrsfs}%
\usepackage[title]{appendix}%
\usepackage{xcolor}%
\usepackage{textcomp}%
\usepackage{manyfoot}%
\usepackage{booktabs}%
\usepackage{caption}
\usepackage{algorithmicx}%
\usepackage{algorithm}
\usepackage{algpseudocode}%
\usepackage{listings}%
\usepackage{cleveref}
\usepackage{algpseudocode}
\usepackage{hyperref}
\usepackage{url}

\newcommand{\beps}{\boldsymbol{\varepsilon}}
\newcommand{\matB}{\underline{\boldsymbol{B}}}
\newcommand{\matL}{\underline{\boldsymbol{\Lambda}}}
\newcommand{\be}{\boldsymbol{e}}
\newcommand{\boldf}{\boldsymbol{f}}
\newcommand{\bg}{\boldsymbol{g}}
\newcommand{\bh}{\boldsymbol{h}}
\newcommand{\bE}{\boldsymbol{E}}
\newcommand{\mean}{\boldsymbol{0}}
\newcommand{\cov}{\underline{\boldsymbol{C}}}  
\newcommand{\matI}{\underline{\boldsymbol{I}}}  
\newcommand{\matD}{\underline{\boldsymbol{D}}}
\newcommand{\bt}{\boldsymbol{t}}
\newcommand{\bu}{\boldsymbol{u}}
\newcommand{\bU}{\boldsymbol{U}}
\newcommand{\bv}{\boldsymbol{v}}
\newcommand{\bV}{\boldsymbol{V}}
\newcommand{\bx}{\boldsymbol{x}}
\newcommand{\bX}{\boldsymbol{X}}
\newcommand{\by}{\boldsymbol{y}}
\newcommand{\bY}{\boldsymbol{Y}}

\theoremstyle{thmstyleone}%
\newtheorem{theorem}{Theorem}
\newtheorem{proposition}[theorem]{Proposition}%

\theoremstyle{thmstyletwo}%
\newtheorem{remark}{Remark}%

\theoremstyle{thmstylethree}%

\begin{document}

\title[Inferring Object Boundaries and their Roughness with Uncertainty Quantification]{Inferring Object Boundaries and their Roughness with Uncertainty Quantification}


\author*[1]{\fnm{Babak} \sur{Maboudi Afkham}}\email{bmaaf@dtu.dk}

\author[1,2]{\fnm{Nicolai Andr\'e Brogaard} \sur{Riis}}\email{nabr@dtu.dk}

\author[1]{\fnm{Yiqiu} \sur{Dong}}\email{yido@dtu.dk}

\author[1]{\fnm{Per Christian} \sur{Hansen}}\email{pcha@dtu.dk}
\equalcont{These authors contributed equally to this work.}

\affil*[1]{\orgdiv{Department of Applied  Mathematics and Computer Science}, \orgname{Technical University of Denmark}, \orgaddress{\city{Kgs.\ Lyngby}, \postcode{2800}, \country{Denmark}}}

\affil[2]{\orgname{Copenhagen Imaging ApS}, \orgaddress{\city{Herlev}, \postcode{2730}, \country{Denmark}}}


\abstract{This work describes a Bayesian framework for reconstructing the boundaries that represent targeted features in an image, as well as the regularity (i.e., roughness vs.\ smoothness) of these boundaries. This regularity often carries crucial information in many inverse problem applications, e.g., for identifying malignant tissues in medical imaging. We represent the boundary as a radial function and characterize the regularity of this function by means of its fractional differentiability. We propose a hierarchical Bayesian formulation which, simultaneously, estimates the function and its regularity, and in addition we quantify the uncertainties in the estimates. Numerical results suggest that the proposed method is a reliable approach for estimating and characterizing object boundaries in imaging applications, as illustrated with examples from X-ray CT and image inpainting. We also show that our method is robust under various noise types, noise levels, and incomplete data.}




\maketitle

\input{1.intro}
\input{2.method}

\section{Numerical Experiments} \label{sec:results}

In this section, we first illustrate the central ideas and methods with a
tutorial example, and then we
evaluate the performance of our method in two difference imaging applications:\
a boundary reconstruction problem in X-ray CT and an inpainting problem. In all numerical tests, the length scale
parameter $\sigma$ in \eqref{eq:lambda} is set as $100$.

\subsection{Data Fitting}\label{sec:datafitting}

To illustrate the concept of roughness we consider time series arising in
electroencephalography (EEG) for analysis of seizures for epileptic patients.
The roughness of the epileptic seizure time series in EEG provides important information
characterizing the type of the seizure \cite{677174,gao2020automatic}.

Hence, we consider the problem of finding the expansion coefficients of a 1D signal
in the Fourier basis from noisy data.
This simple test problem corresponds to the problem \eqref{eq:ip_goal}
with $\mathcal G$ representing the identity operator.

We use the periodic function $\mathcal{V}(x) = x^{3/4}$ defined on the interval $[0,1)$
as the ground truth.
Note that the periodic assumption introduces a discontinuity at $x=0$.
We know from \cite{schwab1998p} that $\mathcal{V}$ belongs to the Hilbert space
$H^{t}_{\text{period}}$ with $t<5/4 - \epsilon$ for all $\epsilon >0$.
Recall the discussion in \Cref{sec:model} that in the limit of the truncation parameter,
i.e., $k\to\infty$, realizations of $\mathcal{V}$ almost surely belong to
$H^{t}_{\text{period}}$ with $t<2s+\nicefrac{1}{2}$.
Therefore, the true roughness parameter $s$ for this signal is at most $0.375$.
The goal of this test problem is to estimate the roughness parameter $s$
together with the signal $\mathcal{V}$ from noisy data.

To discretize the signal $\mathcal{V}$ in order to obtain $\bv$, we use an equally spaced grid $\{(j-1)\Delta x\}_{j=1}^{m}$ with $\Delta x = \ell/m$.
The noise $\bE$ in \eqref{eq:ip_goal} follows an i.i.d.\ Gaussian distribution
with zero mean and standard deviation
$\sigma_{\text{noise}} = r \| \bv\|_2/\| \beps \|_2$ with $\beps$ a realization
of the standard normal distribution and $r$ representing the relative noise level.
Then, the noisy data $\mathbf y$ is obtained following \eqref{eq:ip_deteministic} with $p=m$.


According to the noise model and the problem setting, the likelihood in the posterior formula \eqref{eq:posterior_U} is expressed as
\begin{equation} \label{eq:likelihood_dist}
        \bY | (\bU , S) \sim \mathcal N( \mathcal F (\bu), \sigma^2_{\text{noise}} \matI_p).
\end{equation}
In addition, we assume the prior for $\bU$ is
\begin{equation} \label{eq:prior_dist}
\bU \sim \mathcal N(\mean, \matI_{2k})
\end{equation}
and we assume a non-informative hyper-prior for $S$ in the form of a uniform distribution defined on the interval $[0,10]$.
Then, we can formulate the posterior following \eqref{eq:posterior_U} and compute the samples of $\bV$ through the expression $\bv=\mathcal{F}(\bu)$. In this test problem, we apply the No-U-Turn Sampler (NUTS) 
\cite{JMLR:v15:hoffman14a}
to draw samples from the posterior.
Furthermore, the CUQIpy software package
\cite{CUQIpyPaper}, \cite{cuqipy}
is used to draw $10^4$ samples from the posterior with NUTS\@. We discard $5\cdot 10^3$ samples as burn-in. Python codes demonstrating the experiments in this section is available in \cite{code}.

\subsubsection{Influence of $m$}
\begin{figure}
\centerline{\includegraphics[width=\columnwidth]{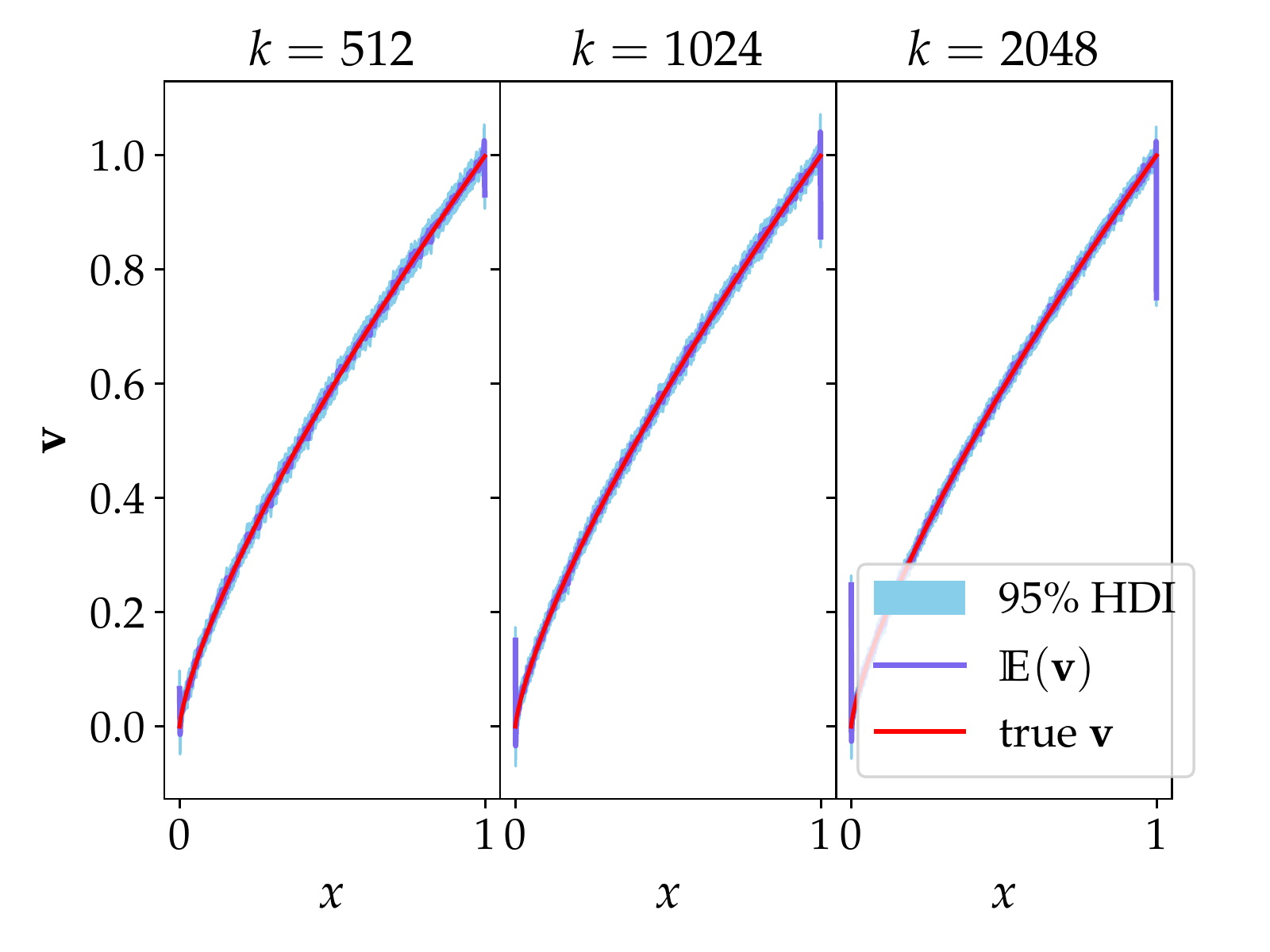} }
\caption{The influence of discretization levels on the reconstructed signals. The plots show the posterior means for the signal $\bv$, and the uncertainties are represented as 95\% HDI.} \label{fig:reg_v_discrete}
\end{figure}

First, we investigate the effect of the discretization level on the reconstruction of $\bv$.
Here, we fix the noise level as $r=0.01$. \Cref{fig:reg_v_discrete} shows the posterior means for the signal $\bv$
with different truncation parameters:\ $k=512$, $1024$ and $2048$. We recall that the discretization level $m$ is set as $m=2k$.
The uncertainty in this estimation is presented in the form of the highest posterior density credibility interval (HDI) \cite{gelman1995bayesian}. We observe that the reconstructed signals are nearly identical to the ground truth and the uncertainties are very low for all discretization levels. Furthermore, it is clear that the reconstruction and the uncertainty are hardly affected by the discretization level. We notice a jump in the reconstructed signals at the boundaries, and this is due to the periodic assumption in the Bayesian prior which results in $\mathcal{V}$ having a discontinuity at the boundaries.

\begin{figure}
\centerline{\includegraphics[width=\columnwidth]{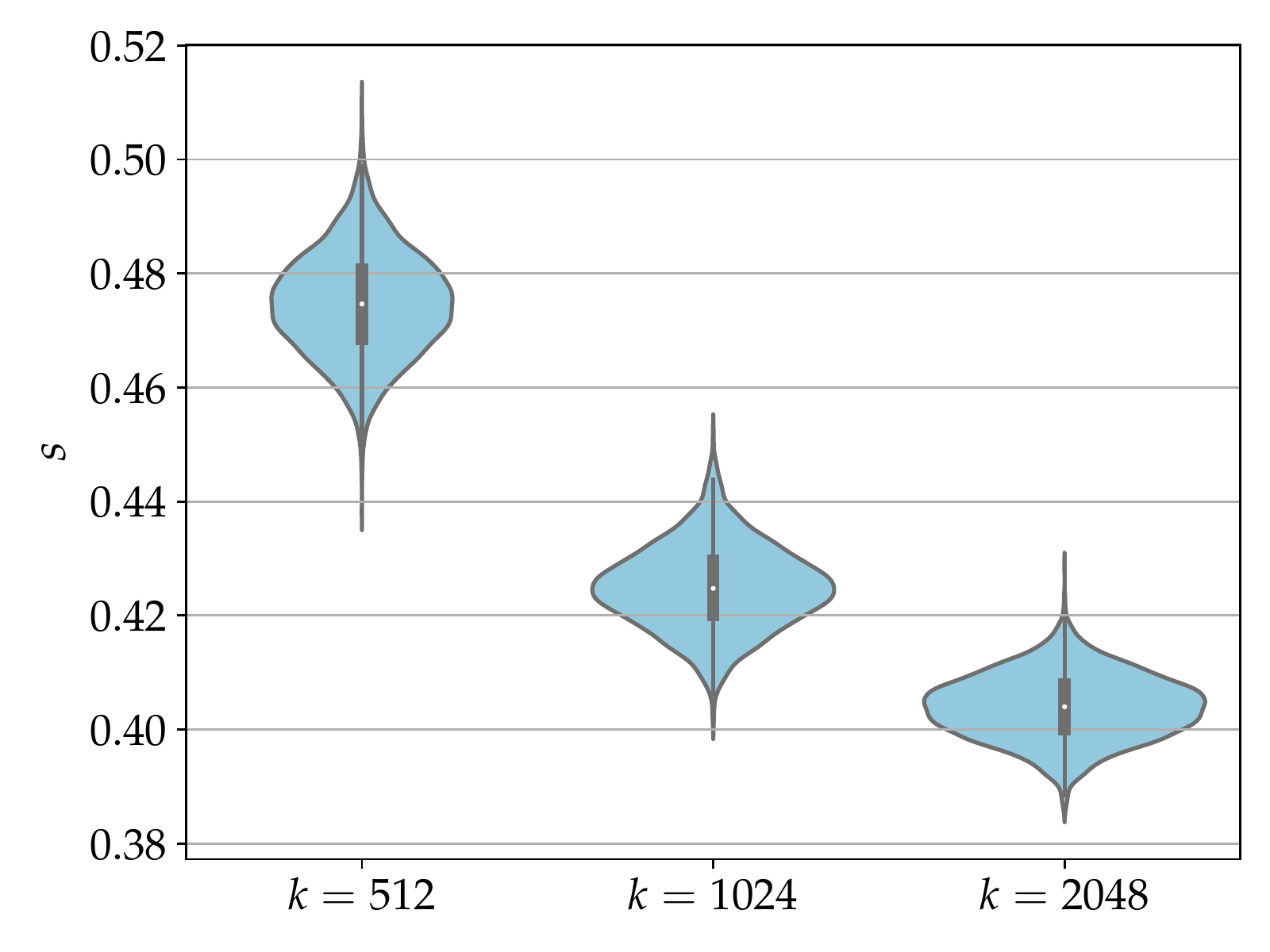} }
\caption{The influence of discretization levels on the roughness parameter $s$. The violin plots show the posterior distributions.} \label{fig:reg_s_discrete}
\end{figure}

In \Cref{fig:reg_s_discrete} we provide violin plots to visualize the uncertainty in estimating $s$. The posterior means for the roughness parameter $s$ are $0.45$, $0.42$ and $0.4$ for $k=512$, $1024$ and $2048$, respectively. Due to the truncation introduced in \eqref{eq:kl_expansion}, we can only expect to estimate an upper bound of the true roughness parameter. Comparing with the fact that the true roughness parameter is at most $0.375$, we can see that the estimated $s$ is a rather tight upper bound for all discretization levels. Furthermore, the uncertainties in $s$ are comparable across all discretization levels.

\subsubsection{Influence of the noise level}

\begin{figure}
\centerline{\includegraphics[width=\columnwidth]{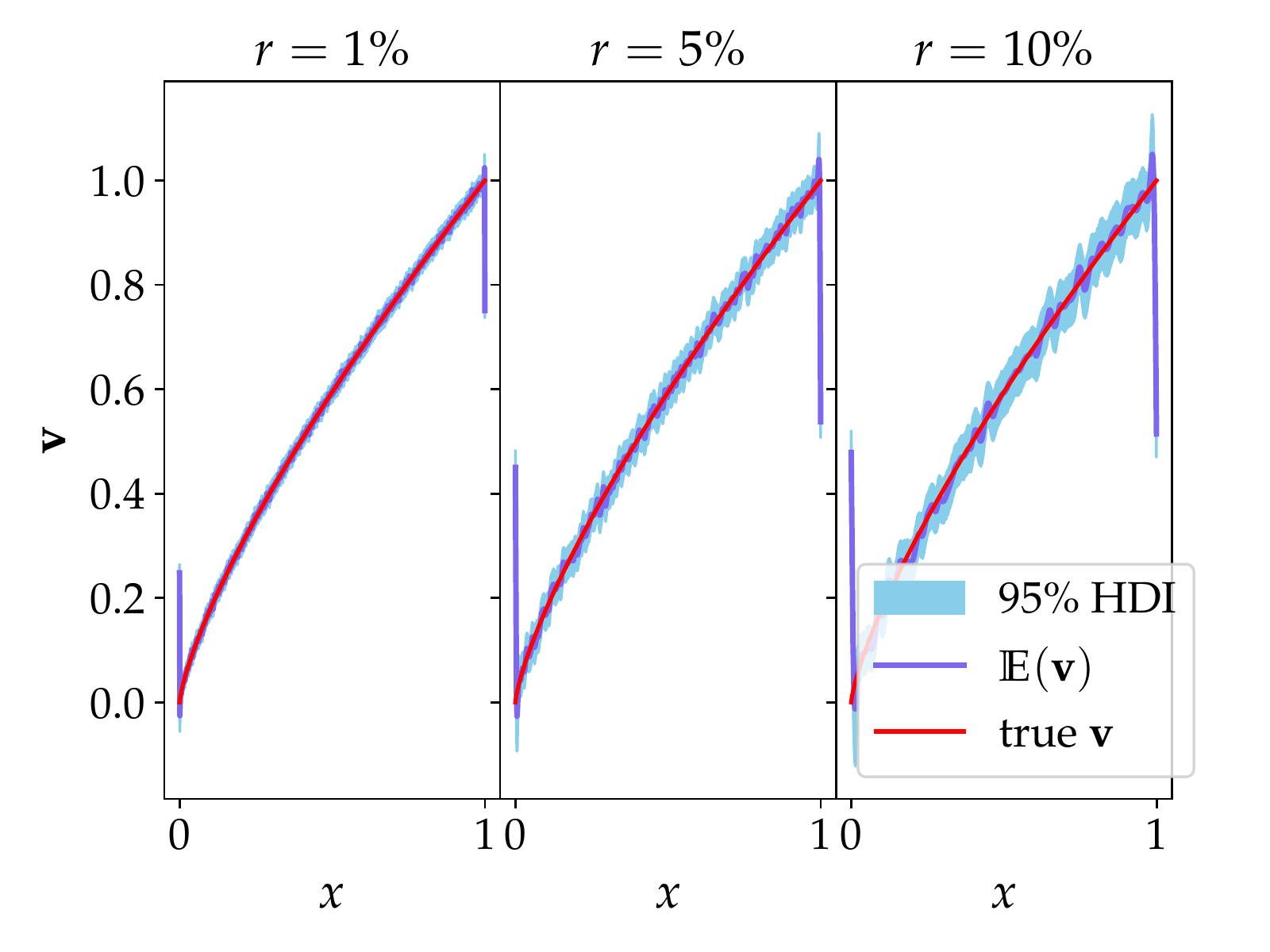} }
\caption{The influence of the noise level $r$ on the posterior statistics for $\bv$. The uncertainty is represented as 95\% HDI.} \label{fig:reg_v_noise}
\end{figure}

We now fix the discretization level to $k=1024$ and investigate the effect of the noise level on the estimation of $\bv$ and $s$. \Cref{fig:reg_v_noise} indicates that the uncertainty in estimating $\bv$ increases with the noise level. However, the reconstructed signal is comparable with the true signal in all cases. 

\begin{figure}
\centerline{\includegraphics[width=\columnwidth]{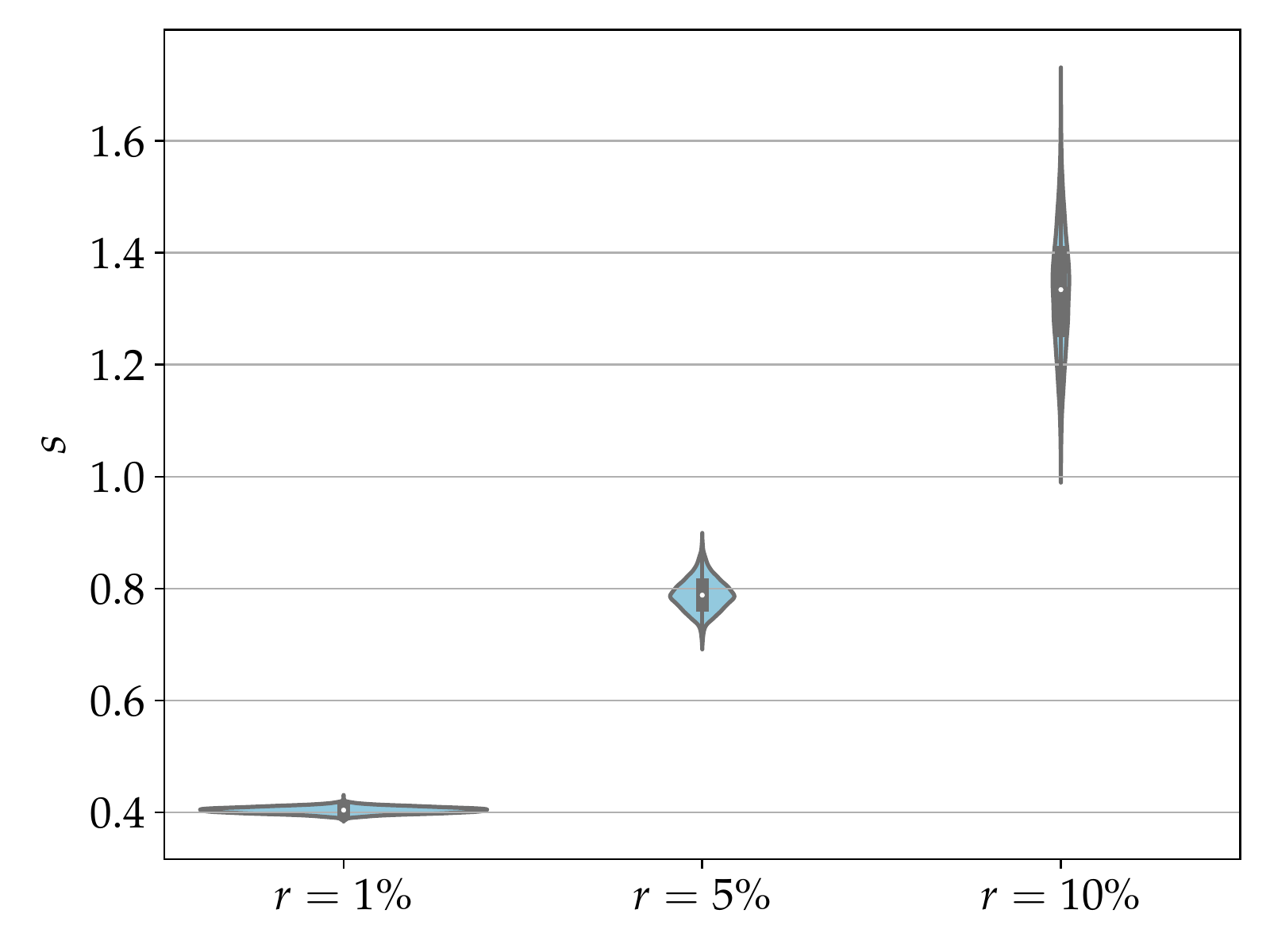} }
\caption{The influence of the noise level $r$ on the posterior statistics for $s$. The violin plots show the posterior distribution $\pi_{S|\bY}$.} \label{fig:re_s_noise}
\end{figure}

The posterior statistics for $s$ are shown in \Cref{fig:re_s_noise}. The posterior means of $s$ are $0.42$, $0.77$ and $1.32$ for the noise levels $r=0.01$, $0.05$ and $0.1$, respectively. It is clear that with the larger noise the estimated upper bound through $s$ becomes looser.
In addition, in \Cref{fig:re_s_noise} we can clearly see that the uncertainty in $s$ increases with the noise level. 

\subsubsection{Comparison with BMS}
Estimating the roughness parameter $s$ in this data-fitting problem can be considered as a model selection problem, where we select one among many candidate models for representing the function $\mathcal{V}$. In this section, we compare our method with a Bayesian model selection (BMS) technique \cite{wasserman2000bayesian,chipman2001practical}. Specifically, we compare with the BMS method for linear regression proposed in \cite{minka2000bayesian}.

To set up the model selection problem, we use $N_{\text{model}}=30$ $s$-candidates that are equally spaced in the interval $[0.75, 2]$, and define the mapping from $\bu$ to $\bv$ as $\mathcal F_j$ according to $s_j$ with $j=1, \cdots, N_{\text{model}}$. BMS is then selecting the mapping $\mathcal{F}$ that gives the largest probability according to the {\it Bayesian evidence}
\begin{equation}\label{eq:evd}
        \pi_{\mathcal F | (\bY = \by) } (\mathcal F_j) \approx \frac{\pi_{\bY| (\bU = \bu^j_{\text{MAP}}, \mathcal F = \mathcal F_j )} (\by)\pi_{(\bU| \mathcal F = \mathcal F_j )} (\bu^j_{\text{MAP}}) }{\sqrt{\det\frac{\matD_j}{2 \pi}}},
\end{equation}
where $\bu^j_{\text{MAP}}$ is the maximum a posteriori (MAP) estimate \cite{kaipio2006statistical} according to $\mathcal F_j$ and $\matD_j=\matL_s(s_j)/\sigma^2_{\text{noise}} + \matI_{2k}$ with $\matL_s$ defined as in \eqref{eq:cov}.

In order to illustrate the sensitivity with respect to the noise, in this experiment we choose rather large noise level $r = 0.1$ and test on two noisy data vectors $\by^1$ and $\by^2$.

\begin{figure}[t]
        \includegraphics[width=\columnwidth]{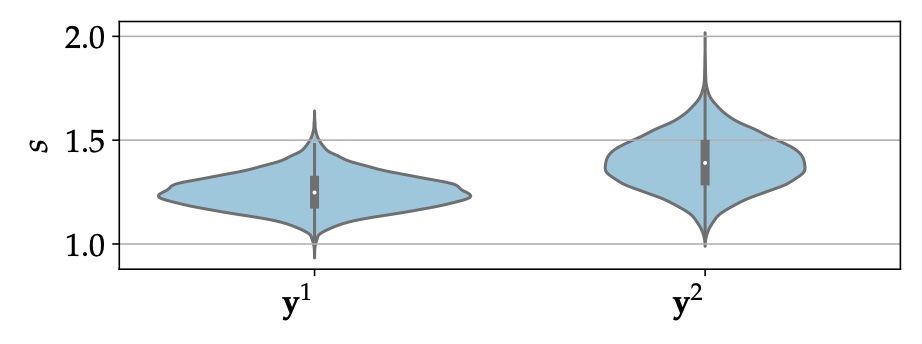} \caption{The posterior distribution $\pi_{S|\bY}$ from our method for $\by^1$ and $\by^2$ with $r=10\%$.} \label{fig:compare_model_selection1}
        \includegraphics[width=\columnwidth]{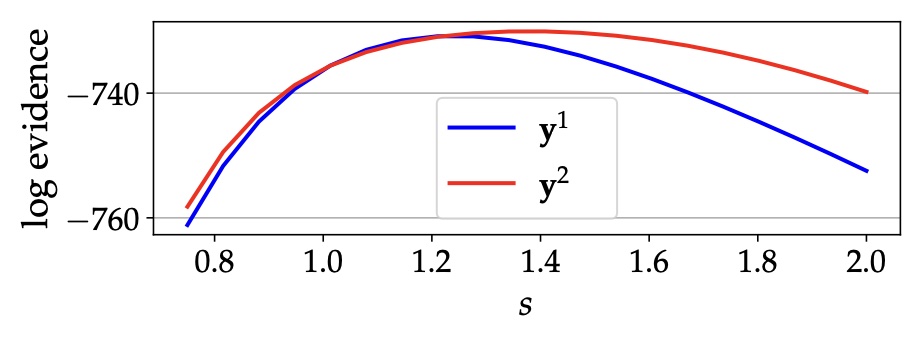} \caption{The logarithm of the evidence \eqref{eq:evd} plotted as a function of $s$ in BMS.} \label{fig:compare_model_selection2}
\end{figure}

In \Cref{fig:compare_model_selection1} and \Cref{fig:compare_model_selection2} we show the performance of both methods. The selected values for the roughness parameter $s$ from BMS, which are according to the peak of the evidence plots, are $1.224$ and $1.397$, respectively. The posterior means of $\pi_{S|\bY}$ from our method are $1.211$ and $1.407$, respectively, and they are very close to the results from BMS. Note that our method provides a distribution of $s$, i.e., $\pi_{S|\bY}$, which can quantify the uncertainty with respect to $\by$. But BMS only shows the influence of $s$ on the evidence, and its output is the value of $s$ according to the peak of the evidence plot. Furthermore, our method does not necessarily require the gradient of the posterior and a discretization of $s$, which are not the case for BMS.



\subsection{Limited-Angle X-Ray Computed Tomography (CT)} \label{sec:ct}

Computed tomography (CT) is a well-established method for reconstructing an image of the cross-section of an object using projection data that represents the intensity loss or attenuation of a beam of X-rays as they pass through the object. According to the Lambert-Beer law \cite[\S 4.2.4]{doi:10.1137/1.9781611976670}, the attenuation of X-rays can be mathematically modelled as the line integrals of the attenuation field in the object, i.e., \textit{Radon transform}. A collection of such line integrals is referred to a \emph{sinogram}. The inverse problem in CT is to reconstruct the attenuation field from a noisy sinogram.

In this example, we examine the performance of our method for identifying the boundary of a single object in the attenuation field together with an estimate of its roughness.
In particular, we are interested in limited-angle imaging configurations, i.e., the projection data are not acquired for all angles around the object.
This could be due to obstacles in the measurement setup or to avoid potentially harmful levels of radiation exposure.
We consider a \emph{parallel-beam geometry} \cite[\S 4.5]{doi:10.1137/1.9781611976670}
with equally spaced projection angles in the interval $[0,\theta_{\text{max}})$,
where $\theta_{\text{max}} \ll \pi$.

We assume that the object has a constant attenuation coefficient on a constant background, and that it follows a \emph{star-shaped} representation as discussed in \Cref{sec:stage}. Here, we assume that the center $\boldsymbol{c}$ is known and recall that $\mathcal A$ denotes the mapping $\mathcal{V}\mapsto \alpha$, $\alpha$ is the attenuation field. We refer the reader to \cite{2107.06607} for inferring the centers and boundaries of multiple objects with a fixed roughness parameter and a non-constant background.

To set up the test problems, we discretize $[0,\theta_{\text{max}})$ to $m$ equally spaced points $\{\iota_j\}_{j=1}^m$ and define $\bv = [\mathcal{V}(\iota_1),\cdots,\mathcal{V}(\iota_m)]^T$ with zero mean.
Then, the forward operator $\mathcal G$ in \eqref{eq:ip_bayes} takes the form $\mathcal G=\mathcal R \circ \mathcal A$, where $\mathcal R$ represents Radon transform. We refer the reader to \cite{doi:10.1137/1.9781611976670} for more detail on X-ray CT\@. In all test cases, we fix $\boldsymbol{c}$ to be the origin, $\alpha^+=2$, $\alpha^-=1$, $r_0 = 0.2$, $b_0 = 0.05$, and $m=256$. We consider Gaussian noise similar to \Cref{sec:datafitting} with $r=0.01$ noise level.

The likelihood and the prior in the posterior \eqref{eq:posterior_U} are defined as
\begin{align}
 \bY | (\bU , S) &\sim \mathcal N(\mathcal{R}\circ\mathcal{A}\circ\mathcal{F} (\bu), \sigma^2_{\text{noise}} \matI_p)\\
\bU &\sim \mathcal N(\mean, \matI_{2k}),
\end{align}
respectively. Furthermore, we consider a uniform distribution on the interval $[0,10]$ for the hyper-parameter $S$. We use the Gibbs sampling method \cite{mcbook} to draw samples from the posterior \eqref{eq:posterior_U}. This method alternatively samples from the conditional distributions
  \begin{equation} \label{eq:conditionals}
  \begin{aligned}
    \pi_{\bU | \bY\!, S }(\bu) & \propto \pi_{\bY|(\bU\!,S)}(\by) \, \pi_{\bU}(\bu) \\[1mm]
    \pi_{S | \bY\!, \bU }(s) &\propto \pi_{\bY|(\bU\!,S)}(\by) \, \pi_{S}(s)
  \end{aligned}
  \end{equation}
by using the Preconditioned Crank–Nicolson MCMC sampler \cite{cotter2013mcmc} and the Metropolis-Hasting sampler \cite{mcbook}, respectively. The sampling algorithm is implemented in the CUQIpy software package. We generate $10^4$ samples from the posterior \eqref{eq:posterior_U} and discard the first $2\cdot 10^3$ samples as burn-in. Python codes demonstrating the experiments in this section is available in \cite{code}.

\subsubsection{Test on a phantom from the prior}
\begin{figure}
\centerline{\includegraphics[width=\columnwidth]{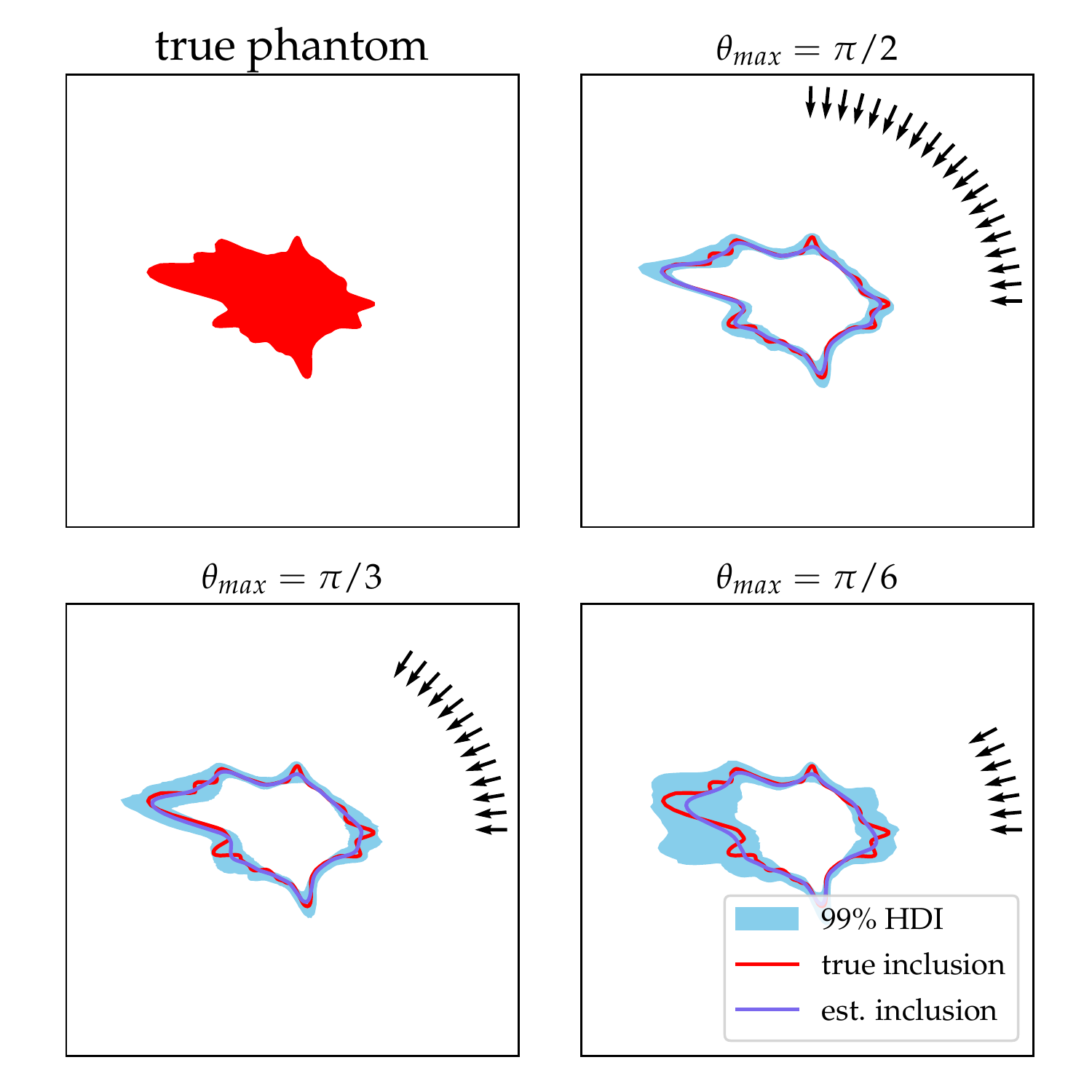} }
\caption{Posterior statistics for the boundary of the object $\bv$ in the limited-angle X-ray CT application; the arrows show the directions of the X-rays.
The ``blob'' shaped phantom (true image) has a boundary that
is chosen from the prior. The red and the dark blue curves indicate the true boundary and the posterior mean, respectively. The light blue region indicates the 99\% HDI.} \label{fig:ct_prior}
\end{figure}

In the first test case, we consider reconstructing a ``blob'' shaped phantom
whose boundary is drawn from the prior distribution, i.e., we draw a sample $(\bu, s)$ from the prior distribution, then compute the vector $\bv$ following \eqref{eq:kl_expansion} and construct a phantom $\alpha$ by using the mapping $\mathcal{A}$. In the top left corner of \Cref{fig:ct_prior}, we show the generated phantom, whose corresponding roughness parameter $s$ is $1.064$. We test our method for $\theta_{\text{max}}=\pi/2$, $\pi/3$ and $\pi/6$ with 384, 256 and 128 projection angles, respectively.
The arrows in \Cref{fig:ct_prior} indicate the angular span of the projections.
The number of detector pixels, i.e. the number of line integrals for each projection angle, is set to 128.

\Cref{fig:ct_prior} shows the posterior means of $\bV$ together with $99\%$ HDI for all three limited-angle cases. We can see that our method provides good estimations of the size and the orientation of the object for all three cases. Furthermore, as the projection interval becomes larger, i.e., $\theta_{\text{max}}$ increases, the accuracy of our method increases and the uncertainty decreases. Even in the extreme limited-angle case with $\theta_{\text{max}}=\pi/6$ our method still can identify some of the spikes in the boundary. In addition, we notice that the true boundary always stays within the uncertainty regions for all $\theta_{\text{max}}$.

Microlocal analysis \cite{Krishnan2015}, \cite[Chapter 8]{doi:10.1137/1.9781611976670}
provides theoretical understanding on which discontinuities in the attenuation field can be stably recovered in CT reconstruction with limited-angle configuration.
This suggests that the uncertainty must be higher at
a part of a boundary if no X-ray is tangent to it. This explains the increase in uncertainty at the boundary curve orthogonal to the projection lines shown in \Cref{fig:ct_prior}.

\begin{figure}
\centerline{\includegraphics[width=\columnwidth]{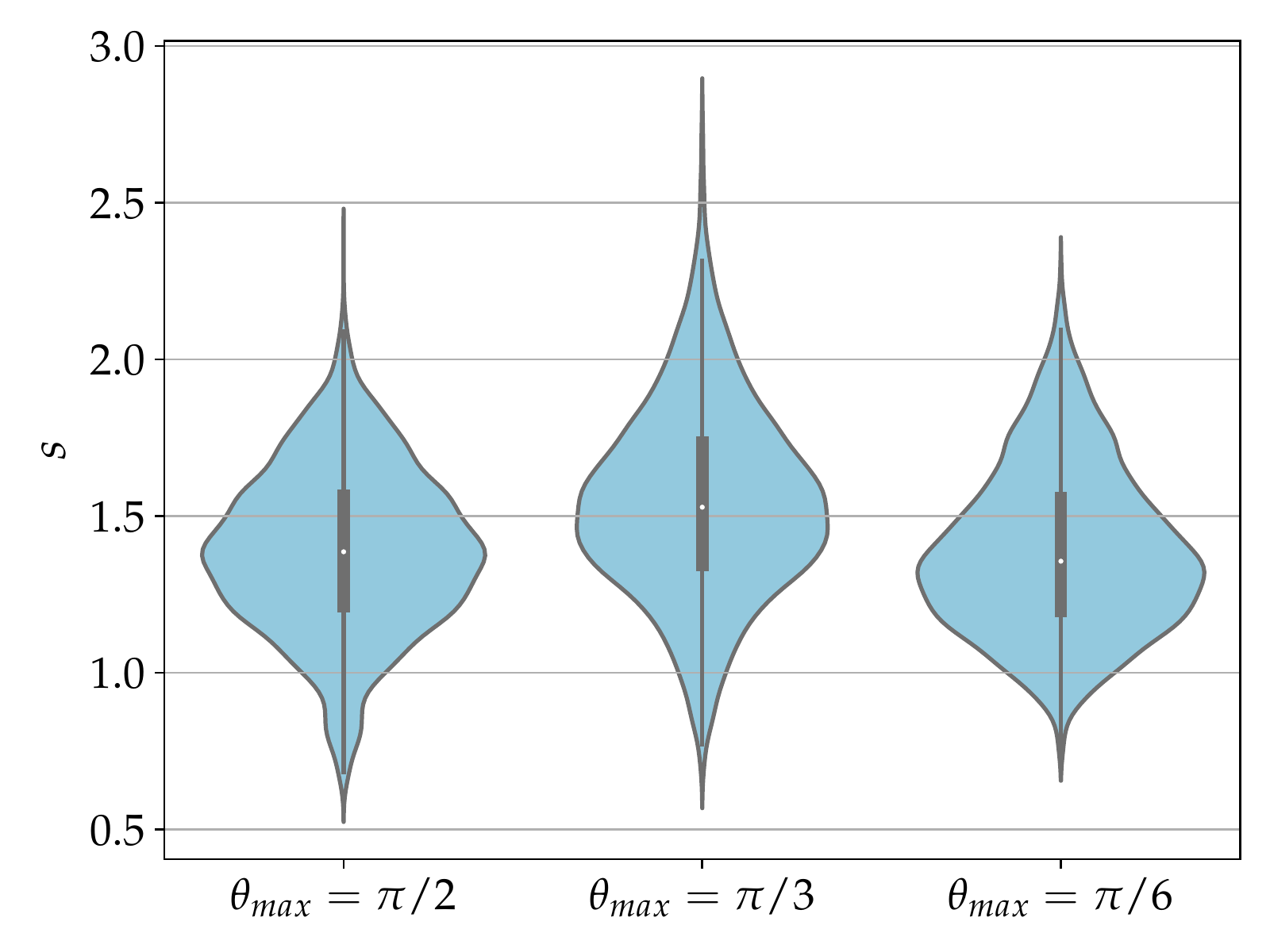} }
\caption{Posterior statistics for $s$ in the case that the phantom is generated from the prior.} \label{fig:ct_prior_s}
\end{figure}

The violin plot in \Cref{fig:ct_prior_s} shows the posterior statistics for the roughness parameter $s$. Comparing with the true value $s=1.064$, the posterior mean of $s$ is a rather tight upper bound. Furthermore, the method provides a comparable estimates of $s$ across all three cases.

\begin{figure}
\centerline{\includegraphics[width=\columnwidth]{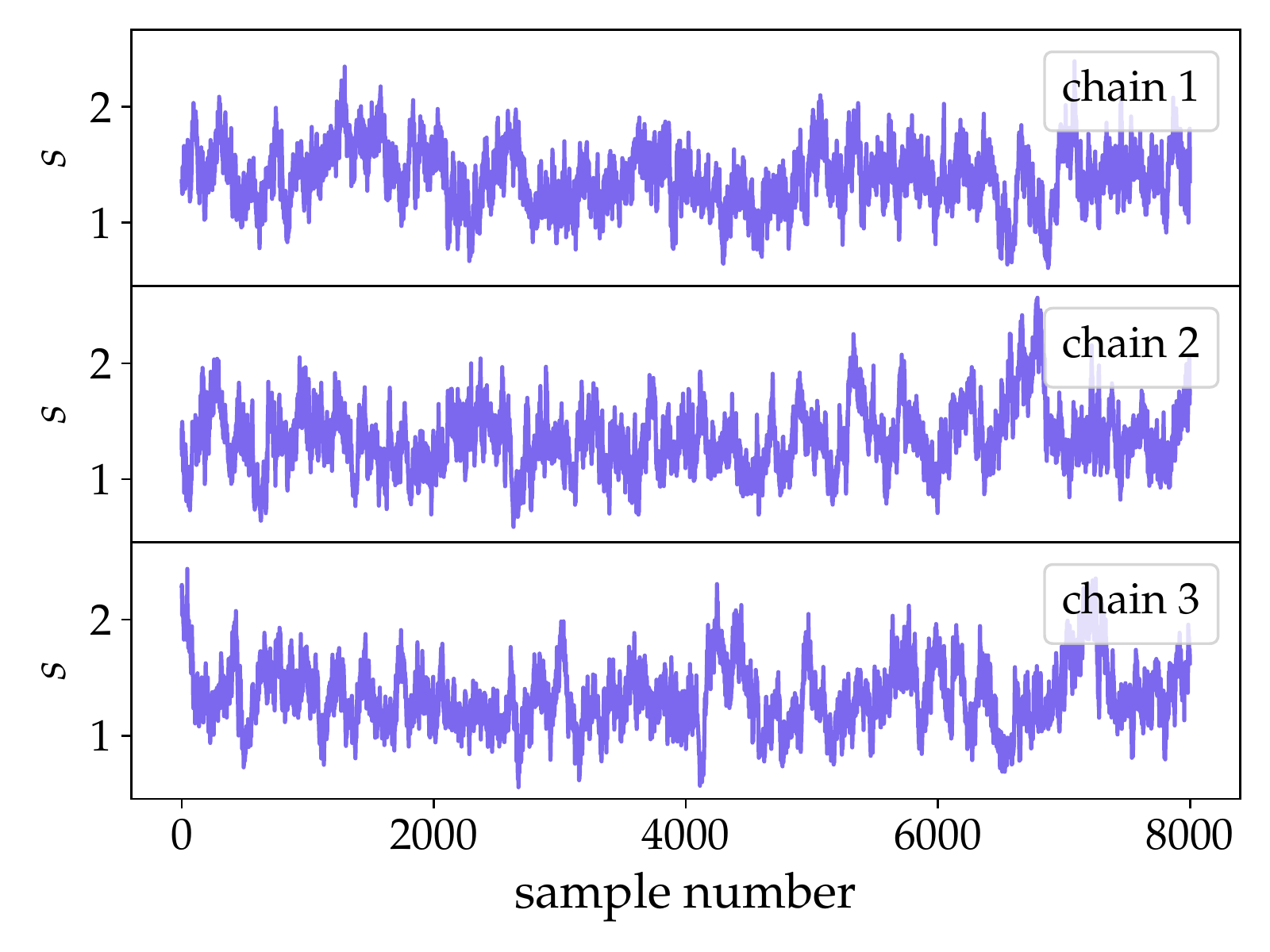} }
\caption{Trace-plots of chains from the Gibbs sampler for $s$.} \label{fig:chains}
\end{figure}

\Cref{fig:chains} shows the trace-plots of the roughness parameter $s$ for 3 chains in the test problem with $\theta_{\text{max}}=\pi/6$. It is obvious that all chains show a rapid mixing, which indicates that a reasonable set of independent sub-samples can be extracted from the chains. We report that the effective sample size (ESS) \cite{mcbook} for the three chains are 87, 132, and 115, respectively. We obtain similar results for the other two cases.

\subsubsection{Test on a gear phantom}
We also test our method on a ``gear'' shaped phantom whose boundary cannot be
represented by the prior.
We consider a phantom with the boundary given by a \emph{gear curve} \cite{weisstein}
\begin{equation} \label{eq:gear_cureve}
        \mathcal{T}_{\text{gear}}(\iota) = r_{\text{gear}}\left(1 + \frac{1}{10} \tanh(10 \sin(n\iota) ) \right), \quad 0\leq \iota < 2\pi.
\end{equation}
Here $r_{\text{gear}}=0.3$ is the mean radius of the gear and $n=10$ indicates the number of gear teeth. The gear phantom $\alpha_{\text{gear}}$ is then constructed following \eqref{eq:phantom} with the center $\boldsymbol{c}$ as the origin.
The gear curve cannot be exactly represented by a finite number of $\sin$-$\cos$ basis functions. In this test, we use the same setup and sampling strategy as in the previous test problem except taking $m = 512$.

\begin{figure}
\centerline{\includegraphics[width=\columnwidth]{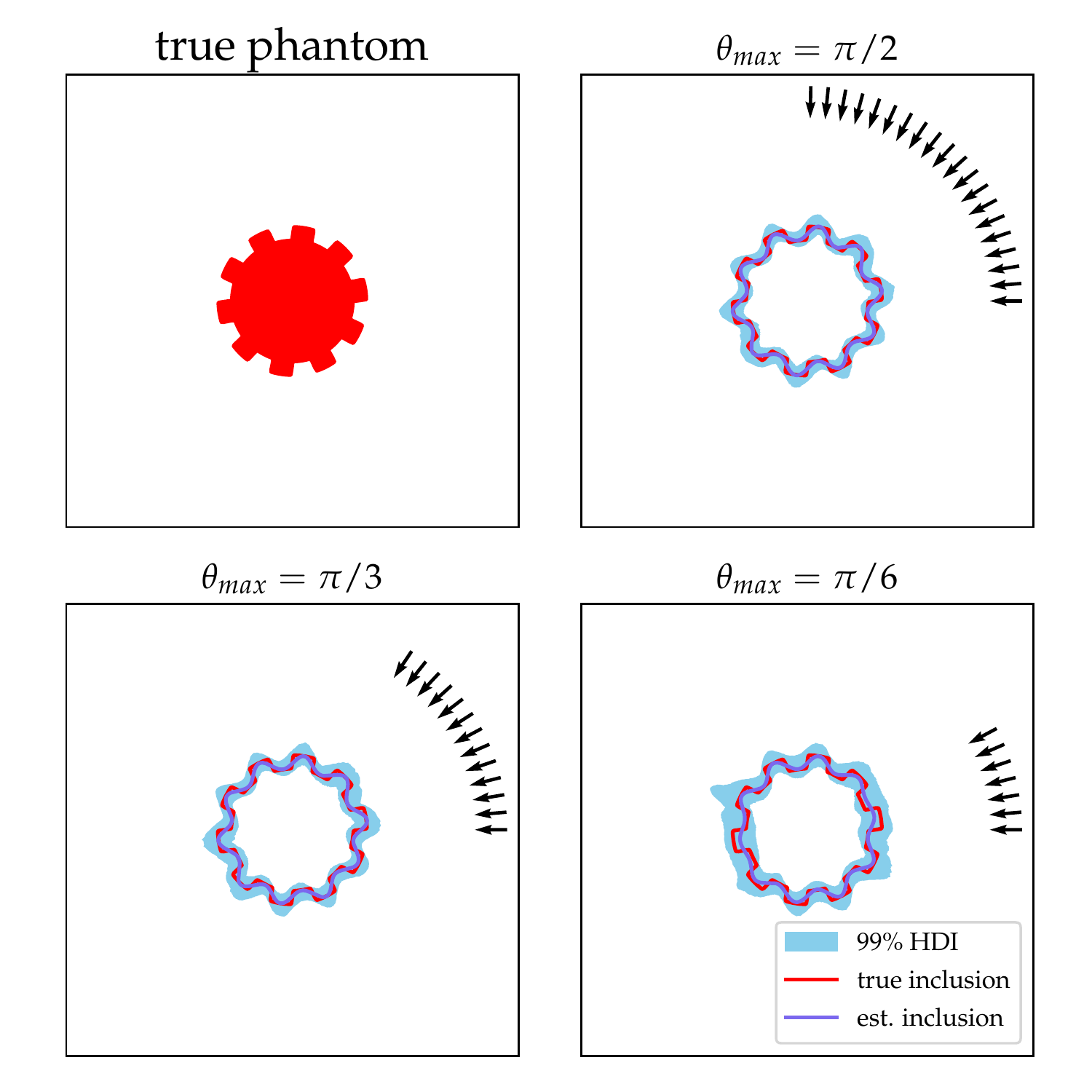} }
\caption{Posterior statistics for $\bv$ from the phantom that is constructed from a gear curve. The red and the dark blue curves indicate the true boundary and the posterior mean, respectively. The light blue region indicates the 99\% HDI.} \label{fig:ct_gear}
\end{figure}

\begin{figure}
\centerline{\includegraphics[width=\columnwidth]{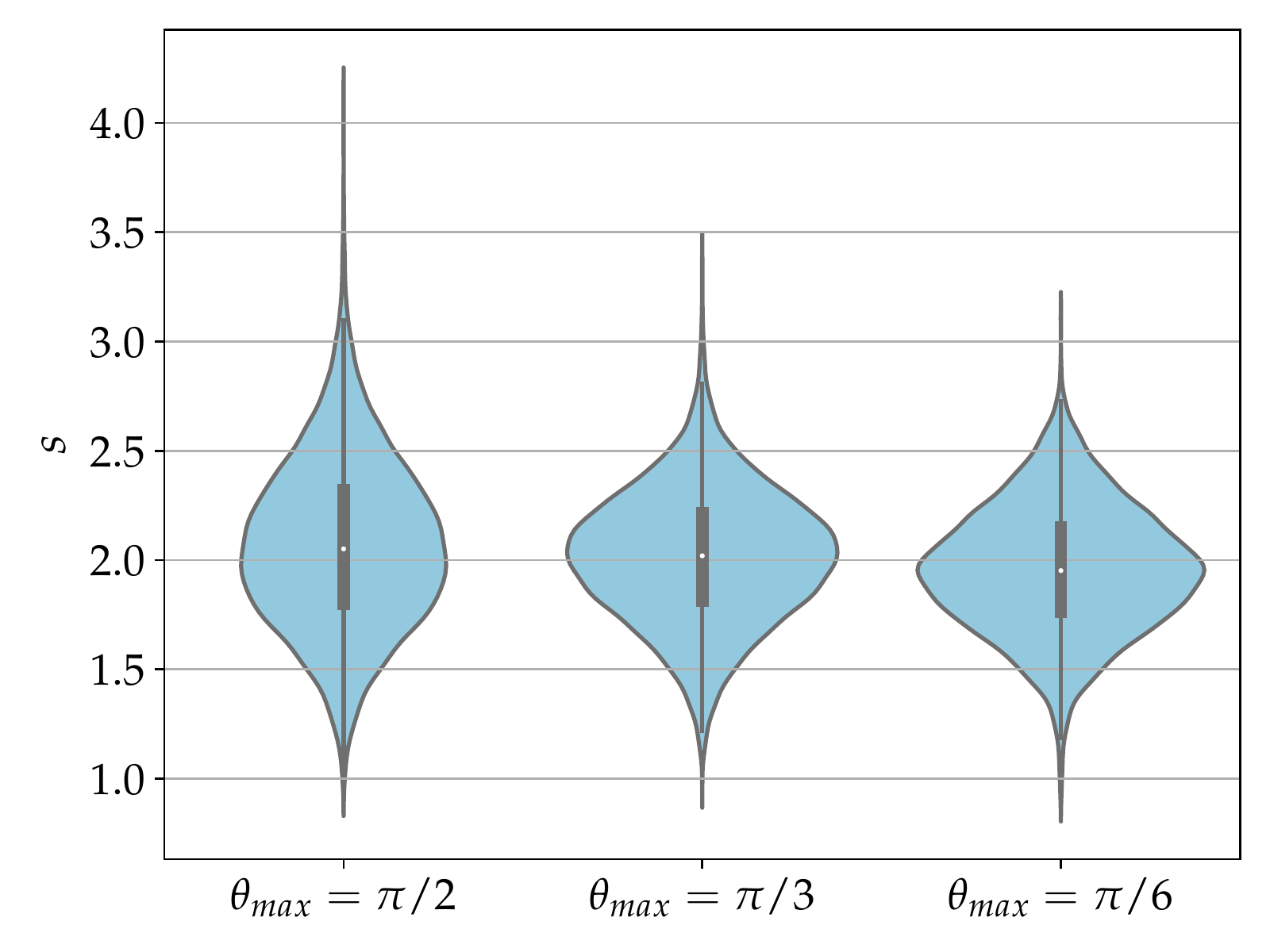} }
\caption{Posterior statistics for the roughness parameter $s$ for the phantom generated from the gear curve.} \label{fig:ct_gear_s}
\end{figure}

\Cref{fig:ct_gear} shows the posterior statistics for the gear boundary $\bv$ under the same three limited-angle cases. We can see that the location of the gear teeth are recognized in all cases. Further, the sharp corners of the gear teeth are all rounded in the reconstruction. The reason is that sharp corners cannot be represented by a finite number of the sin-cos basis. Comparing the results for different $\theta_{\text{max}}$ we also notice that the uncertainty increases as $\theta_{\text{max}}$ decreases, and the part of the boundary that no projection is tangent to it has larger uncertainty. In addition, we show the posterior distribution of the roughness parameter $s$ in \Cref{fig:ct_gear_s}. It is obvious that the posterior means in all three cases are comparable.

\subsubsection{Comparison with multi-step method}
We compare our method with the multi-step method, i.e., first reconstructing the 2D attenuation field then segmenting the object and detecting the boundary. We use the filtered back projection (FBP) method \cite{doi:10.1137/1.9781611976670} for the reconstruction. To detect the boundary, we apply thresholding-based method to segment the object and then use the Sobel filter to obtain the boundary \cite{castleman1996digital}. We test the cases with $\theta_{\text{max}} = \pi/2$ as well as $\pi/3$ and the noise level $r=0.01$. We omit the comparison for the case with $\theta_{\text{max}} = \pi/6$, since the multi-step method cannot provide meaningful results.

%
%

\begin{figure}[t]
\begin{center}
\begin{minipage}[t]{3.5cm}
        \includegraphics[totalheight=3.3cm]{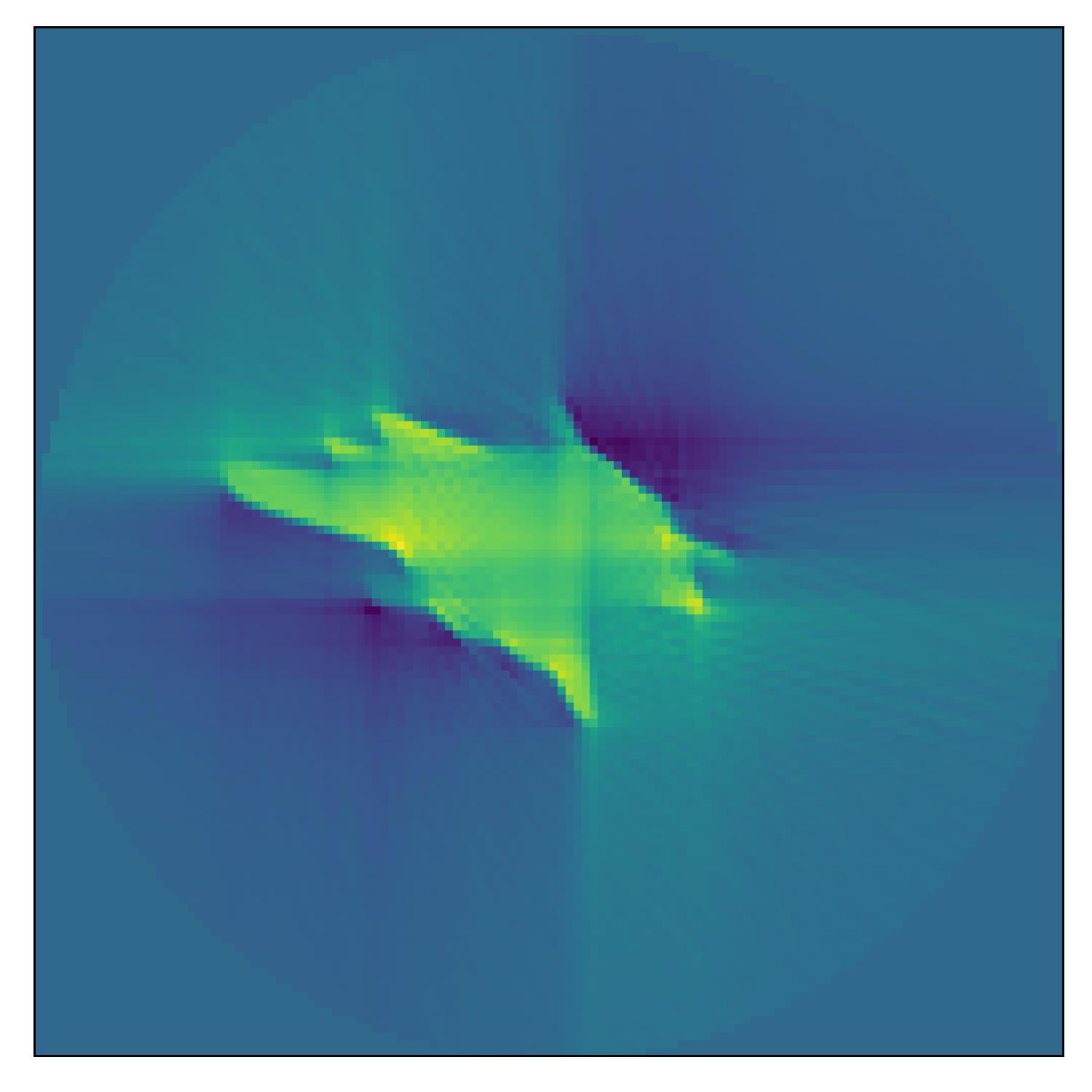} \\
        \centering{(a)}
\end{minipage}
\begin{minipage}[t]{3.5cm}
        \includegraphics[totalheight=3.3cm]{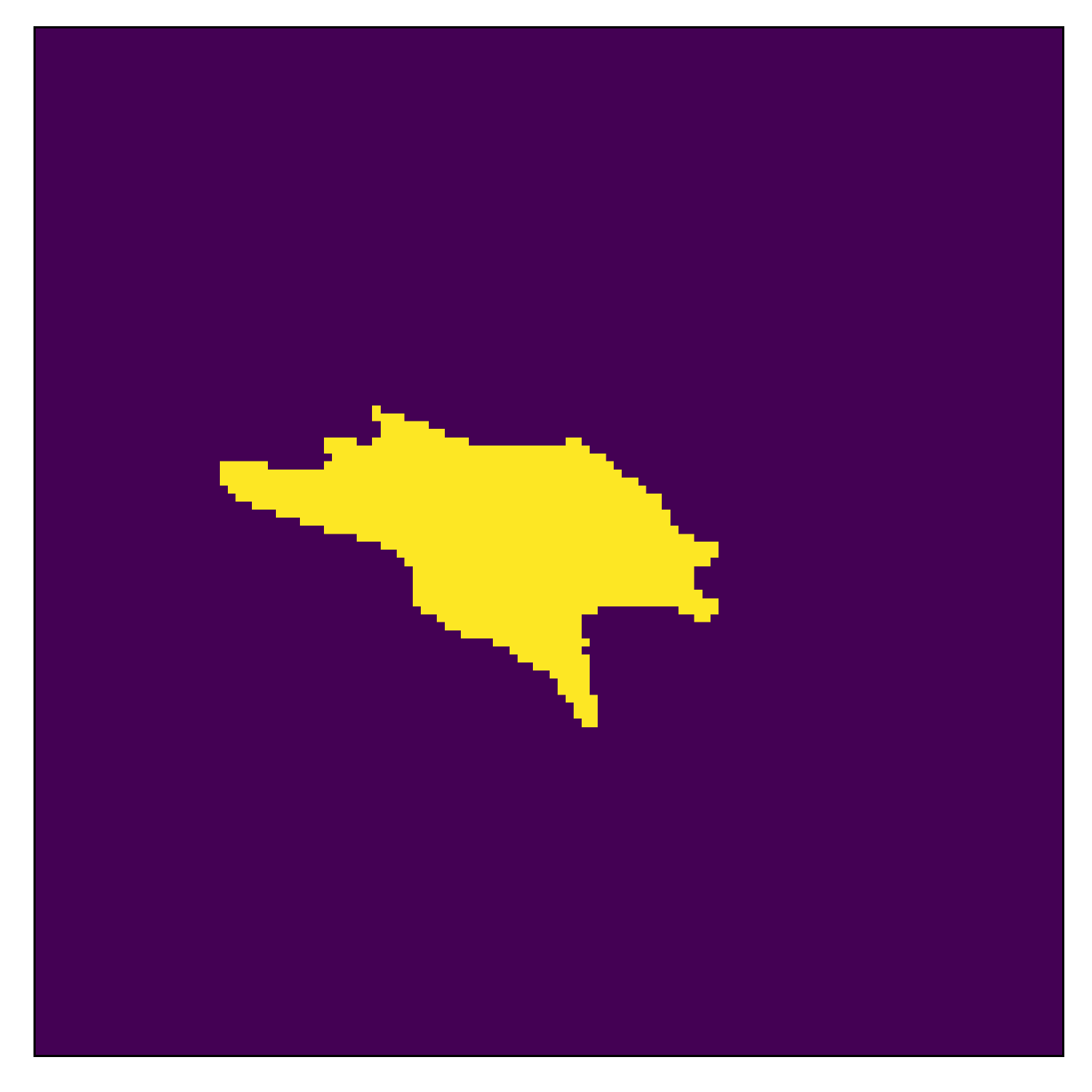} \\
        \centering{(b)}
\end{minipage}\\
\begin{minipage}[t]{3.5cm}
        \includegraphics[totalheight=3.3cm]{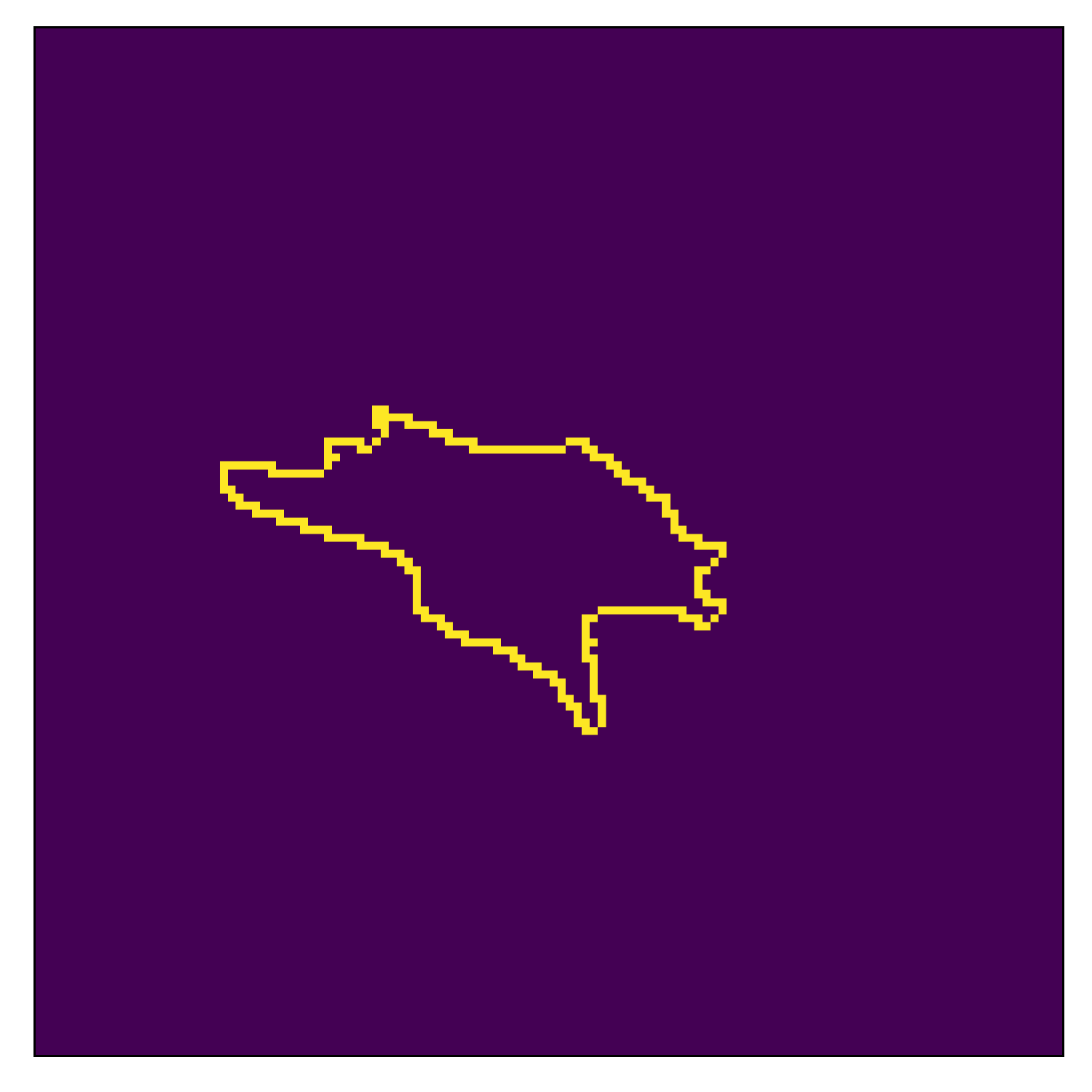} \\
        \centering{(c)}
\end{minipage}
\begin{minipage}[t]{3.5cm}
        \includegraphics[totalheight=3.3cm]{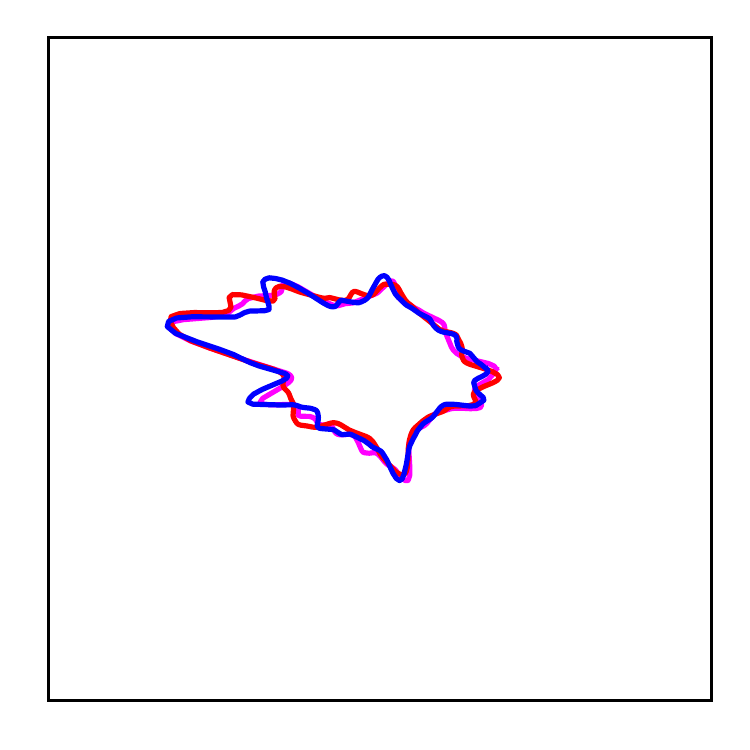} \\
        \centering{(d)}
\end{minipage}
\end{center}
        \caption{Comparison of the multi-step method and our method for the case $\theta_{\text{max}} = \pi/2$. The results from the multi-step method: (a) FBP reconstruction; (b) image segmentation; (c) boundary detection. The results from our method: (d) 3 samples from $\pi_{\bV|\bY}$.} \label{fig:fbp_90}
  \end{figure}

  \begin{figure}[t]
  \begin{center}
\begin{minipage}[t]{3.5cm}
        \includegraphics[totalheight=3.3cm]{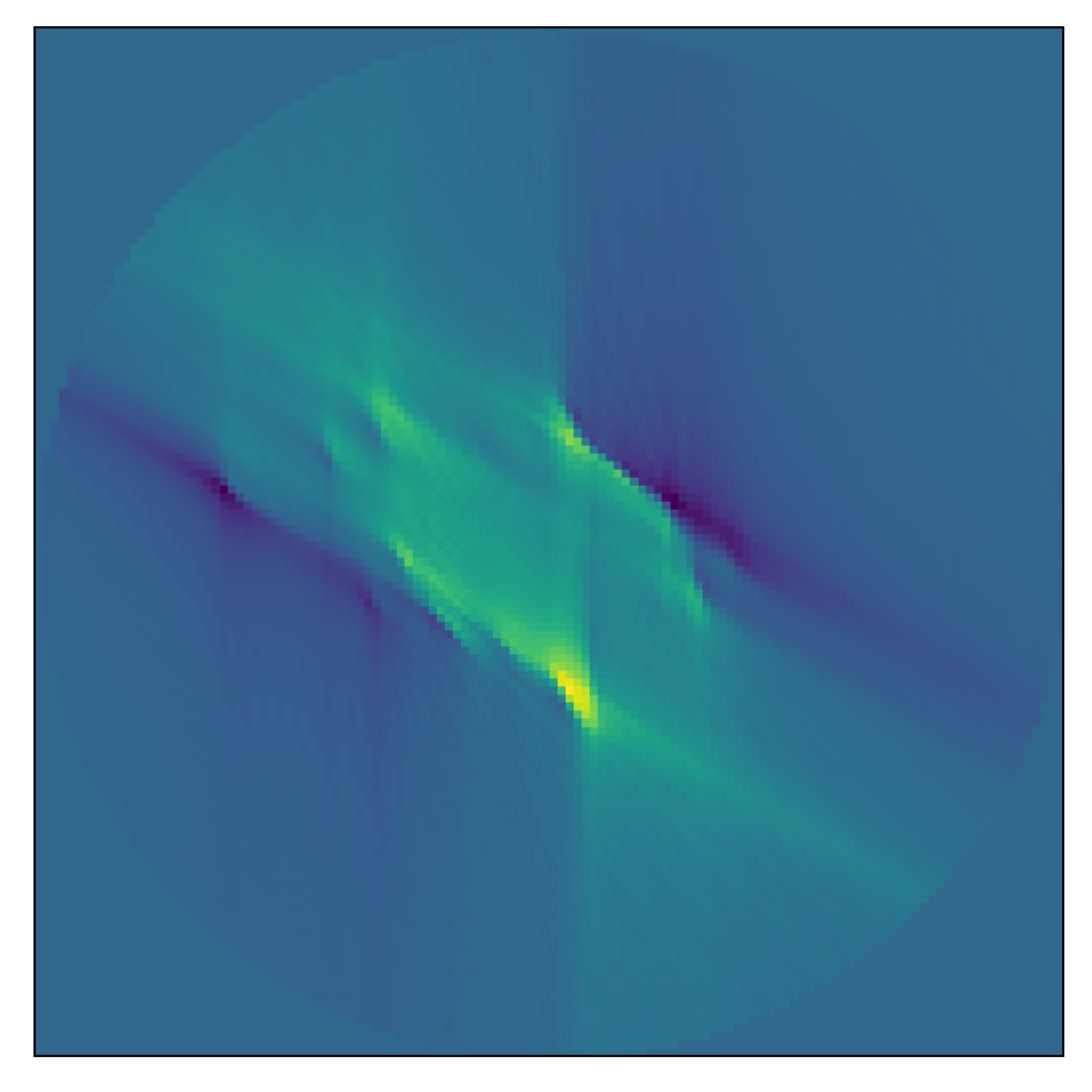} \\
        \centering{(a)}
\end{minipage}
\begin{minipage}[t]{3.5cm}
        \includegraphics[totalheight=3.3cm]{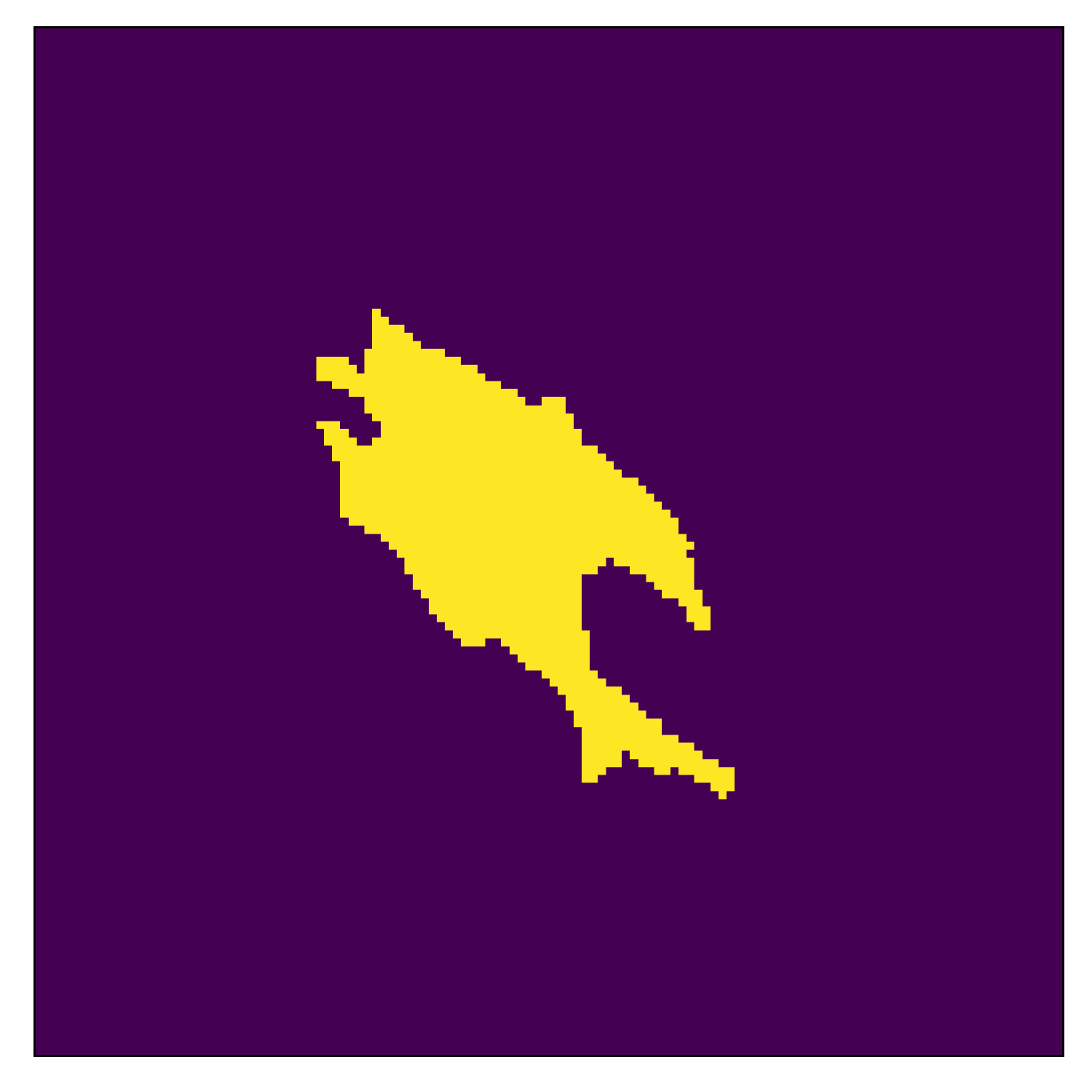} \\
        \centering{(b)}
\end{minipage}\\
\begin{minipage}[t]{3.5cm}
        \includegraphics[totalheight=3.3cm]{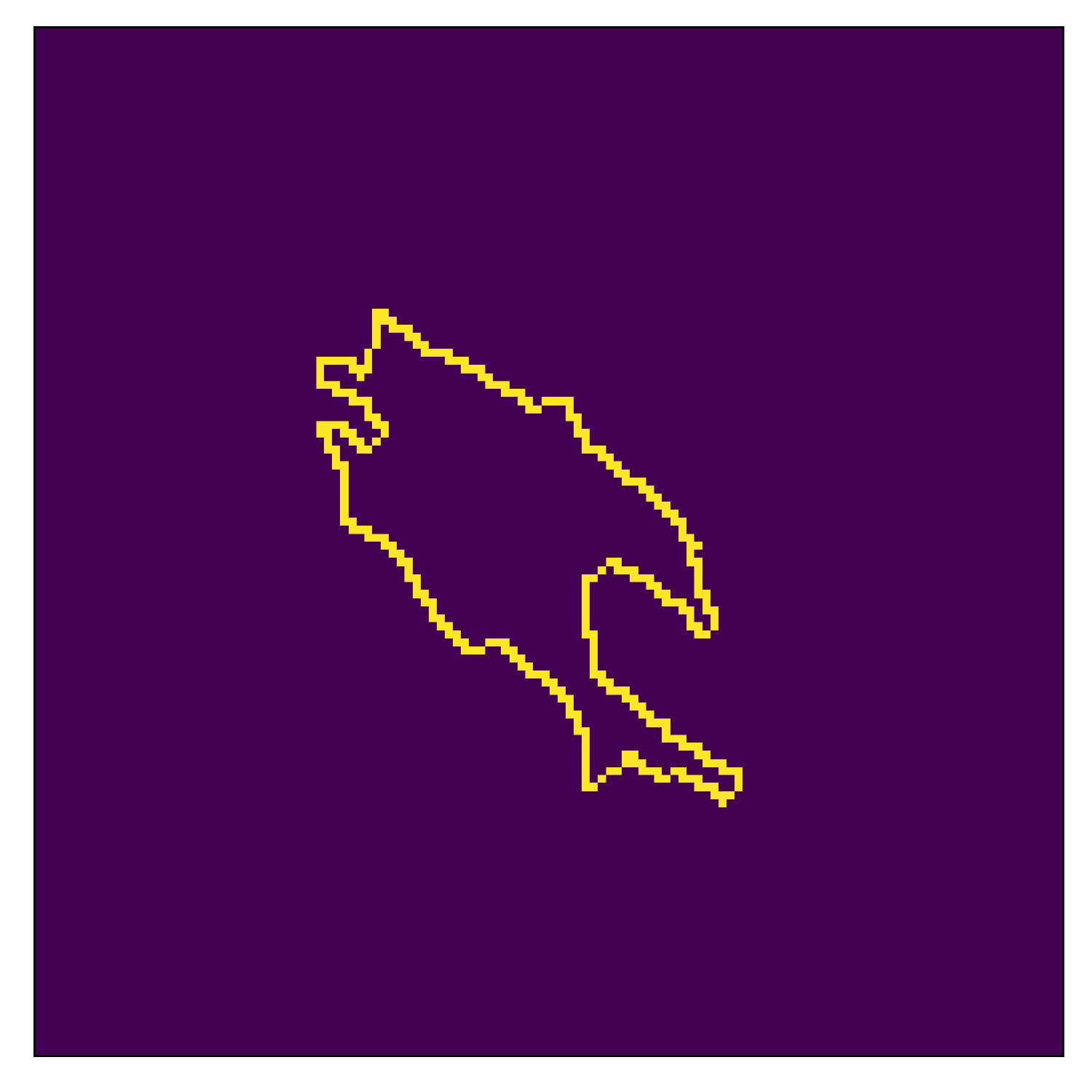} \\
        \centering{(c)}
\end{minipage}
\begin{minipage}[t]{3.5cm}
        \includegraphics[totalheight=3.3cm]{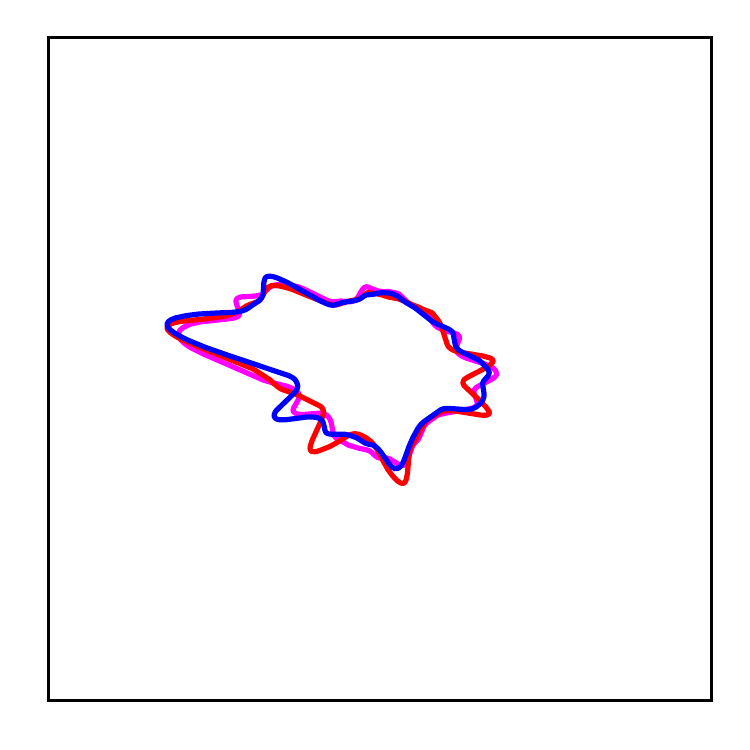} \\
        \centering{(d)}
\end{minipage}
\end{center}
        \caption{Comparison of the multi-step method and our method for the case $\theta_{\text{max}} = \pi/3$. The results from the multi-step method: (a) FBP reconstruction; (b) image segmentation; (c) boundary detection. The results from our method: (d) 3 samples from $\pi_{\bV|\bY}$.}  \label{fig:fbp_60}
\end{figure}

In \Cref{fig:fbp_90} and \Cref{fig:fbp_60} we show the results from the multi-step method in all steps, and compare with 3 samples from our method. Note that the posterior means and $99\%$ HDI from our method can be found in \Cref{fig:ct_prior}. For the multi-step method, we can clearly see that the artifacts in the reconstruction due to limited-angle scenarios continue affecting the segmentation and boundary detection. We note that there are many advanced CT reconstruction methods and edge detection methods, see, e.g., \cite{doi:10.1137/1.9781611976670, 1596801}. However, without incorporating more specific prior knowledge on the boundary, the multi-step methods would face the same challenges and provide similar results for such limited-angle scenarios. In addition, the multi-step method can only provide a single result without any quantification of uncertainties. We show 3 samples from the posterior distribution of our method in \Cref{fig:fbp_90}(d) and \Cref{fig:fbp_60}(d), and the difference among samples in fact highlights the larger uncertainty regions. Further, it's obvious that our method handles the challenges from limited-angle scenarios much better.


\subsection{Image Inpainting}

To conclude our experiments we test our method on an image inpainting application. The goal is to recover the missing features from damaged and degraded images \cite{schonlieb2015partial}. We mainly investigate how our method performs under different noise models.


We consider inferring the boundary of a unique object from an incomplete noisy image together with the roughness of the boundary. We assume that the intensities of the object and the background are constants and the boundary of the object has a star-shaped description given in \eqref{eq:phantom} with $\alpha^+ = 1$ and $\alpha^- = 0$. The other parameters in the star-shaped description are chosen identical to the setup in \Cref{sec:ct}. We set $m=256$ for the discretization of the boundary and $256\times 256$ as the image size. The mapping from the boundary to the image is denoted by $\mathcal A: \mathcal{V} \mapsto \alpha$. The forward operator $\mathcal G$ in \eqref{eq:ip_bayes} takes the form $\mathcal G = \mathcal L \circ \mathcal A$, where $\mathcal{L}$ is the data loss operator.
Note that for numerical stability, the boundary of the object is slightly smoothed, and this smoothness is taken into account in the operator $\mathcal{L}$.

\begin{figure}
\centerline{\includegraphics[width=\columnwidth]{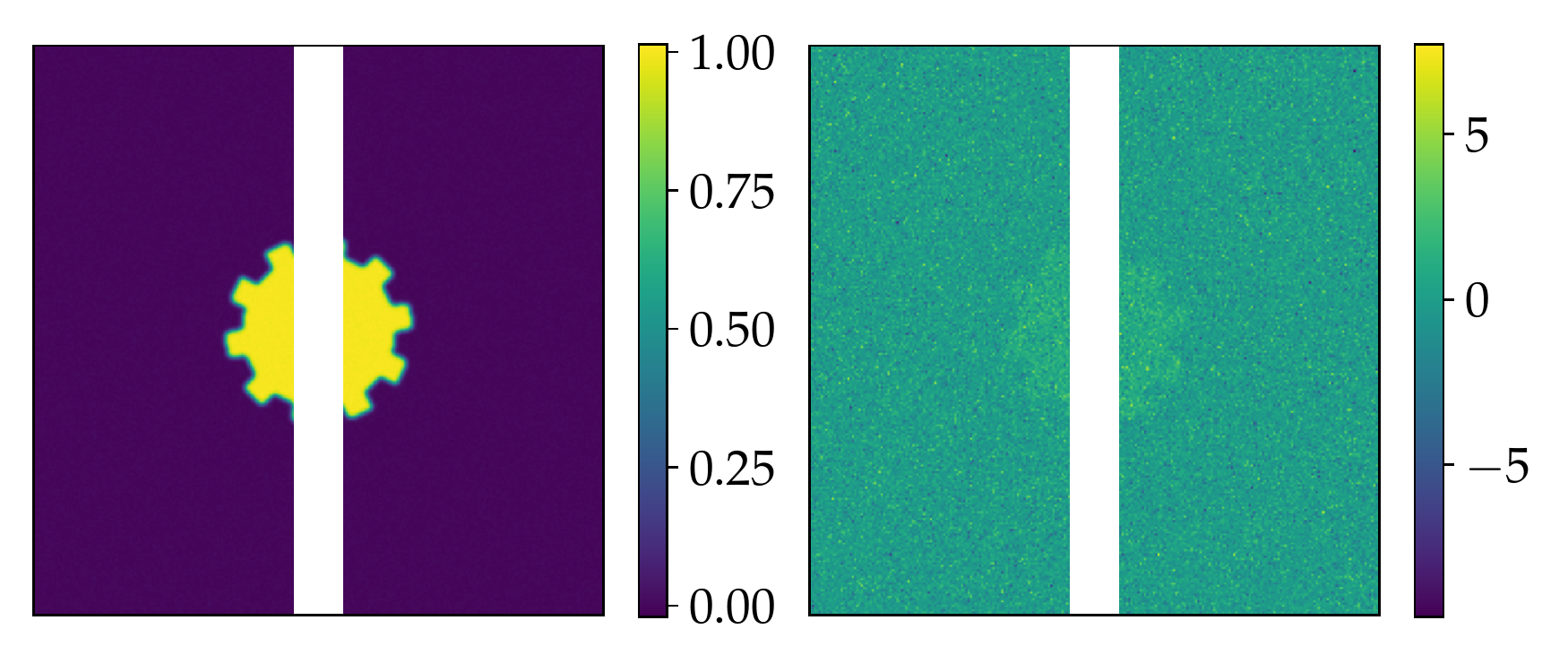} }
\caption{The degraded images corrupted by Gaussian noise (left) and Laplace noise (right), respectively. The white stripes indicate the missing regions.} \label{fig:incomplete_obs}
\end{figure}

To investigate the performance of our method under different noise models, we consider Gaussian noise with the noise level $r = 0.02$ and Laplace noise. The latter is another type of additive noise that often models impulse noise \cite{L1TV}. In this case, $\bE$ follows a multivariate Laplace distribution with mean $\mean$ and
 covariance matrix $\sigma^2_{\text{noise}} \matI_m$, and $\sigma^2_{\text{noise}}= 0.8$.
We show both degraded images in \Cref{fig:incomplete_obs}. The white stripes indicate the regions where the information of the pixels is lost.


We use the same sampling method as in \Cref{sec:ct}. For the Gaussian noise case, we draw $10^4$ samples with $10^3$ samples as burn-in, and the ESS value is 131. For the Laplace noise case, we need more samples due to high-tail behavior of Laplace distribution, and we draw $2\cdot 10^4$ samples with $10^3$ samples as burn-in. Furthermore, we notice that the posterior in this case turns out to be multi-modal,
so ESS is not an effective diagnostic estimate for the sampling procedure. We leave chain diagnostics for such a distribution to future studies.

\begin{figure}
\centerline{\includegraphics[width=\columnwidth]{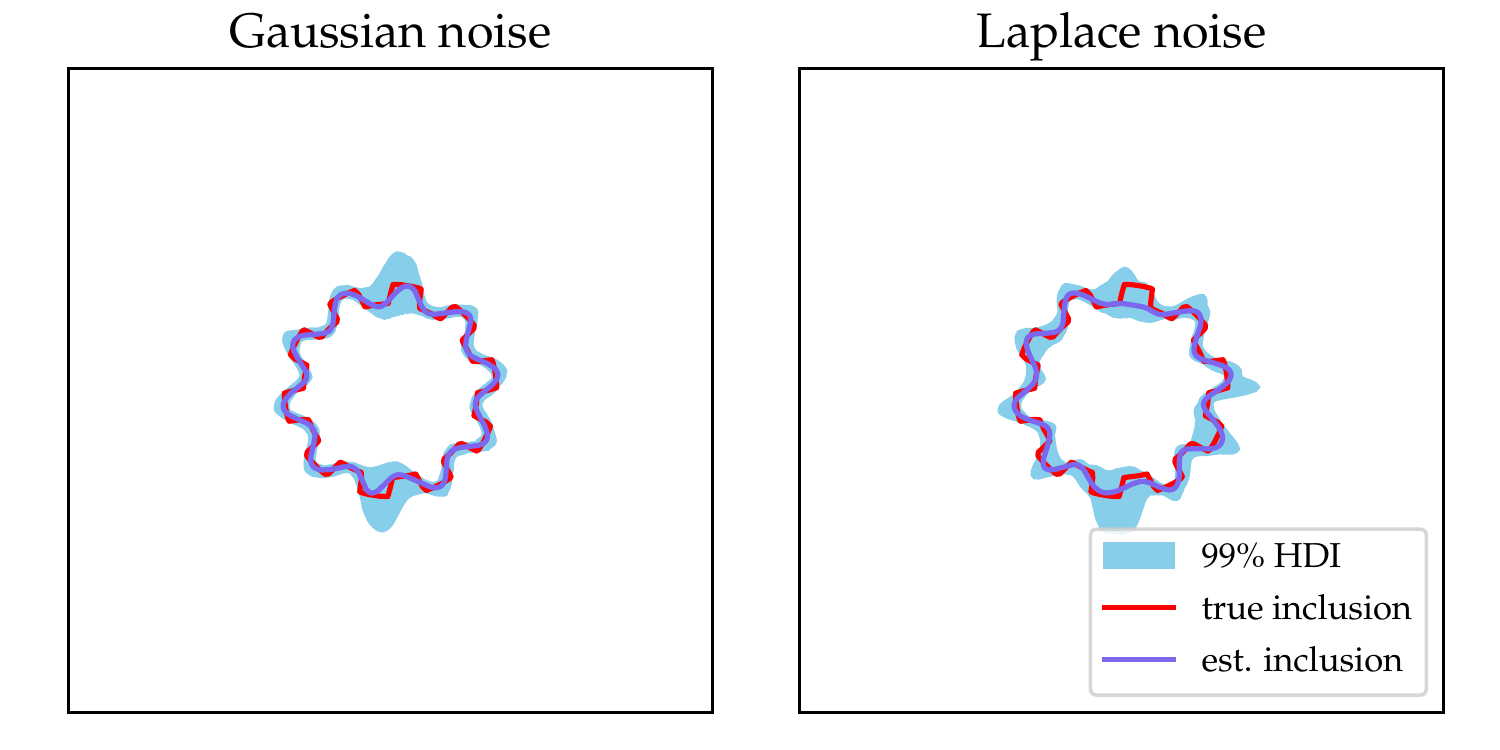} }
\caption{Restored boundaries together with the uncertainties in image inpainting application. The red and the dark blue curves indicate the true boundary of the gear and the posterior means, respectively. The light blue regions indicate the 99\% HDI.} \label{fig:incomplete_gear}
\end{figure}

\Cref{fig:incomplete_gear} shows the restored boundaries together with their uncertainty under both noise models. The restored boundaries are obtained from the posterior mean, and the uncertainties of the restorations are presented with 99\% HDI.

For the Gaussian case, all teeth of the gear are identified in the restored boundary, and larger uncertainty appears at the location of the missing pixels. Further, the uncertainty band contains the true boundary.

For the Laplace case, all teeth except the one located at the top of the gear in the missing region are correctly identified.
The top tooth is not identified in the posterior mean, and the uncertainty band
suggests that there is indeed a large amount of uncertainty in this estimate.
Furthermore, we notice that the true boundary curve is contained in the uncertainty band. Due to the multi-modality of the posterior, the posterior mean is not the best representation of the true boundary. The maximum a posteriori (MAP) estimate may be a better choice. The computation of the MAP estimate for such a posterior is left as future research.

\begin{figure}
\centerline{\includegraphics[width=\columnwidth]{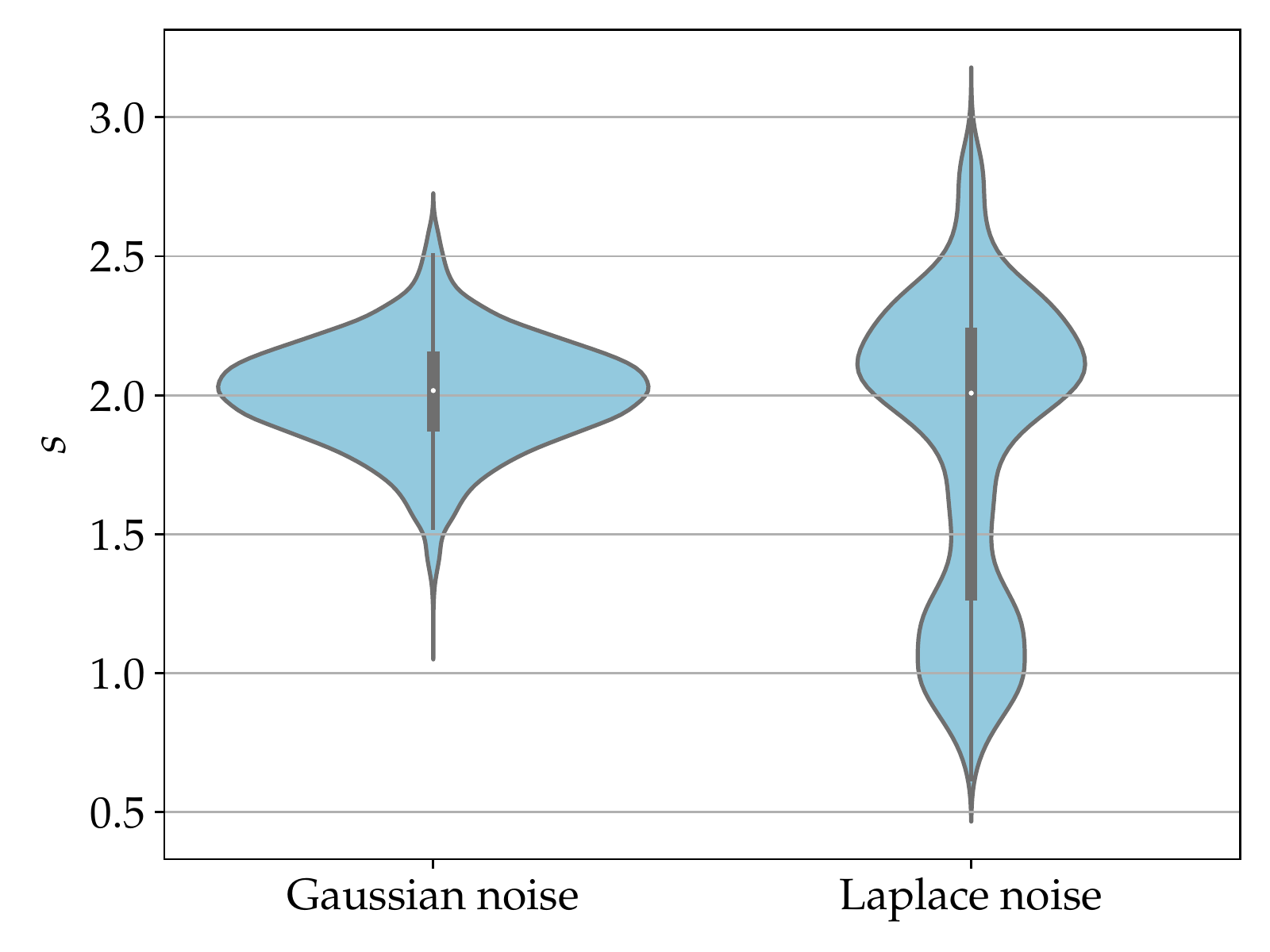} }
\caption{Posterior statistics for the roughness parameter $s$ in image inpainting application.} \label{fig:incomplete_s}
\end{figure}

The posterior statistics for the roughness parameter $s$ are presented in \Cref{fig:incomplete_s}. For the Gaussian case, the posterior mean for $s$ is 2.01, which is comparable with the estimation shown in \Cref{fig:ct_gear_s}. For the Laplace case, we notice that the distribution is multi-modal, which is a strong suggestion that the posterior for $(\bV, S)$ is multi-modal as well.
Although the posterior mean of $s$ is estimated to be $1.83$, we notice that the dominant mode in the posterior distribution of $s$ is around $2.2$.
This is fairly consistent with the estimate for the regularity of the gear curve from the case with Gaussian noise, similar to the finding in \Cref{sec:ct}.

\section{Conclusions} \label{sec:conclusions}

This paper presents a novel Bayesian framework for reconstructing functions with
uncertain regularity, which has important applications in various inverse problems,
and we applied the framework to the computation of object boundaries in imaging applications.
Our method simultaneously estimates both the function and its regularity by characterizing
its roughness by means of fractional differentiability.
We then present a hierarchical Bayesian formulation of this problem.
By exploring the posterior distribution, we can estimate the unknown function
and quantify uncertainties in the estimation.
We also demonstrate how our method can be implemented efficiently by means
of the FFT and IFFT algorithms.

Numerical results indicate that this Bayesian framework is a promising and effective tool for solving a broad range of inverse problems where the regularity of functions carries critical information. Furthermore, our experiments showcase the excellent performance of the method in quantifying the uncertainties for inverse problems with various noise types, noise levels, and incomplete measurement.

Although we focus on inferring 1D features with uncertain roughness, we can naturally extend our work to estimating features in higher dimensions. Furthermore, for problems that contain uncertainties regarding other parameters of the Whittle-Mat\'ern covariance, e.g., the length scale $\sigma$, we can construct a similar posterior and explore it using a Gibbs sampler with the structure presented in this paper. As an interesting direction of future research, we propose generalizing the findings of this paper to other priors with hierarchical structure, e.g., Besov priors \cite{Dashti2012}.

\section*{Acknowledgement}

We thank Professor Alessandro Foi, from Tampere University, Finland, for comments to an earlier version of the manuscript
that helped to improve the presentation. This work was funded by a Villum Investigator grant (no. 25893) from the Villum Foundation.

\bibliography{sn-bibliography}

\end{document}

%% file: 1.intro.tex

\section{Introduction}

A very important task in imaging science is to reconstruct geometric features of an image, such as the shape or boundary of an object within the image.
In this work we focus on computation of the boundary of an object, and the goal is to
determine the boundary as well as its roughness.
Moreover, we describe how to compute the uncertainty of these features due to
errors in the data.

A classical approach to identify the boundary is to first reconstruct the image
from the data and then
apply a post-processing step, e.g., image segmentation.
This has two drawbacks:\
there is a waste on computations,
and errors in the
reconstructed image may propagate to the segmentation and thus the computed boundary.
The latter is pronounced in the case of X-ray computed tomography (CT)
with few projections or with beam hardening,
where reconstruction artifacts are common \cite{doi:10.1137/1.9781611976670}.
Therefore, it is natural to avoid error propagation by combining the two steps
as done in, e.g., \cite{ALouis2011, Ramlau2007, PCH2017} and this is also our approach.

There are recent works that combine such steps for inverse problems involving integral equations, e.g., in Fourier imaging and magnetic resonance imaging (MRI), X-ray computed tomography (CT), and deconvolution applications. As a non-exhaustive list of recent works see \cite{ye2006asymptotic,ye2003cramer,4303091,4531632}. These methods infer geometric aspects of objects, e.g., shape, size and orientation, and quantify the uncertainties in such geometric parameters. However, these methods still lack the estimation and uncertainty quantification (UQ) of the roughness of the boundaries, with respect to noisy data.

Knowledge of boundary roughness is crucial in many application. For example, in the context of X-ray CT, the roughness of the boundaries of tumors
is essential for differentiating between benign and malignant tissues
\cite{limkin2019complexity,nakasu1999preoperative,sanghani2019evaluation,baba2012tumour}. Confidence in estimating such features helps with an accurate diagnosis and designing an effective treatment plan.
However, understanding the uncertainties in estimating the roughness using quantitative methods
is not well-established.

The present work builds on our previous method from \cite{2107.06607}, where we focused solely on
CT reconstruction problem and computing the boundary and its uncertainty. The main contribution of \cite{2107.06607} is a theoretical study
in an infinite-dimensional setting.
We now consider the general inverse problem
in
a finite-dimensional setting, and focus on the numerical computation of the boundary as well as its roughness directly from the data. Moreover, we perform UQ on both features.
\emph{The main contribution of this work is a Bayesian framework for a generic inverse problem}, and a fast algorithm for inferring the
boundary and its roughness as well as their uncertainties.

The key ingredient is a hierarchical Bayesian approach to estimating the roughness of a
scalar function, based on a Whittle-Mat\'ern covariance prior and the
level of fractional differentiablity.
We also present an efficient FFT-based computational method for implementing this.

We illustrate the central ideas and the method with
a data-fitting problem related to electroencephalography.
We then evaluate its performance through two different imaging problems:\
a boundary reconstruction problem in X-ray CT, and an image inpainting problem.
Our numerical experiments suggest that the proposed method is a reliable approach
for inferring object boundaries.
We show that our method has a robust performance under various levels of measurement noise,
different noise models, as well as for challenging scenarios with incomplete data.

Throughout the paper, we use the following notation.
Probability distributions are denoted $\pi_{\Box}$ where the subscript
is used to characterize the specific distribution.
We indicate a real value with a lower case letter $x$,
a real-valued vector with a bold letter $\bx$, a real-valued random variable
(associated with $x$) with a capital letter $X$, and a real-valued multivariate random
variable (associated with $\bx$) with a bold capital letter $\bX$.
Moreover, a constant and deterministic matrix is denoted by an underlined boldface
capital letter $\underline{\bX}$.

Our paper is organized as follows.
We introduce the Bayesian inverse problem and a hierarchical formulation for inferring
boundaries with uncertain roughness in Sections \ref{sec:stage} and \ref{sec:model}.
Fast computational algorithms 
are developed in \Cref{sec:method}.
In \Cref{sec:results} we illustrate the performance of our method,
and we present concluding remarks in \Cref{sec:conclusions}.

\section{A Goal-Oriented Bayesian Approach}
\label{sec:stage}

We assume that the boundary of an object is characterized by a real-valued function
on a continuum, $\mathcal{V}:\mathbb R \to \mathbb R$,
and we let the vector $\bv \in \mathbb R^m$ be a discrete representation
of $\mathcal{V}$.
Then the
forward model for the
inverse problem of reconstructing the boundary takes the form
  \begin{equation} \label{eq:ip_deteministic}
    \by = \mathcal G(\bv) + \beps,
  \end{equation}
where $\by \in \mathbb R^p$ is the measured data and $\beps \in \mathbb R^p$ is the noise.
Furthermore, $\mathcal{G}:\mathbb R^m \to \mathbb R^p$ is the \emph{forward operator} that describes the relation between $\bv$ and $\by$.
This way of defining a 2D image from its 1D boundary
can significantly reduce the number of unknown variables in the inverse problem,
which is commonly referred to as a \emph{goal-oriented} strategy for this specific inverse problem.

As an example of such representation of a boundary, we present the problem of finding tumor boundaries in the context of X-ray CT. We assume that tumors have a constant attenuation coefficient on a constant background, and that it follows a \emph{star-shaped} representation \cite{2107.06607,Dunlop2016}. This means that the boundary of the object is expressed by a center $\boldsymbol{c}=(c_x, c_y)^T$ and a radial function $\mathcal{T}(\iota) = r_0 + b_0 \exp(\mathcal{V}(\iota))$ with $\iota\in[0, 2\pi)$, where $r_0>0$ is the minimum radius, $b_0>0$ is a scaling parameter, and $\mathcal{V}$ is a periodic function on $\iota\in[0,2\pi)$.
With a given boundary, the attenuation field $\alpha(x,y)$ is defined by
  \begin{equation} \label{eq:phantom}
    \alpha(x,y) =
        \left\{
        \begin{aligned}
                \alpha^+, \quad & \sqrt{ (x-c_x)^2 + (y-c_y)^2 } < \mathcal{T}(\iota^*), \\
                \alpha^-, \quad & \text{otherwise} ,
        \end{aligned}
        \right.
  \end{equation}
where $\iota^* = \arctan( (y-c_y)/(x-c_x) )$, and $\alpha^+,\alpha^-\geq0$ are the attenuation coefficients of the object and the background, respectively. Let $\mathcal A$ denotes the mapping $\mathcal{V}\mapsto \alpha$. The forward operator $\mathcal G$ in \eqref{eq:ip_deteministic} takes the form $\mathcal G=\mathcal R \circ \mathcal A$, where $\mathcal R$ represents a discrete Radon transform, constructing an X-ray attenuation measurement, i.e., a sinogram $\by$. In the rest of this article we consider inverse problems of reconstructing boundaries of objects in a piece-wise constant field, such as reconstructing tumor boundaries in X-ray CT discussed above.

The Bayesian approach to inverse problems is a promising way to quantify
uncertainties in the solution \cite{kaipio2006statistical}.
In this approach we model all quantities of interest in \eqref{eq:ip_deteministic}
as random variables, each with their own probability distribution.
Thus, we introduce the random variables $\bY$, $\bV$, and $\bE$ and consider
the stochastic model for the inverse problem formulated as
  \begin{equation} \label{eq:ip_bayes}
    \bY = \mathcal{G}(\bV) + \bE .
  \end{equation}
The Bayesian solution, called the \emph{posterior}, is then the conditional distribution
of the boundary given some data, and it is denoted by $\pi_{\bV | \bY }(\bv)$.
The characteristics of the posterior, e.g., its moments, are used to interpret the
uncertainty associated with the solution of the inverse problem.

The Bayesian approach
has been applied to a wide range of inverse problems in imaging, such as
X-ray CT \cite{2107.06607},
electrical impedance tomography \cite{Matthew_2016,Dunlop2016},
digital holography \cite{carpio2022seeing},
and hydraulic tomography \cite{xu2020preconditioned,alghamdi2020bayesian,alghamdi2021bayesian2}.

In order to characterize the roughness of the boundary and quantify its uncertainty,
we use the \emph{level of fractional differentiability} \cite{DINEZZA2012521}.
The advantage of such characterization is that the behavior of the boundary
is independent of the discretization.
The \emph{Whittle-Mat\'ern priors} \cite{1930-8337_2014_2_561} are a class of
Gaussian priors that provides control over fractional differentiability.
Hence we assume that $\bV$ follows a Whittle-Mat\'ern prior in the form
$\bV \sim \mathcal N (\mean, \cov_s)$, i.e.,
$\bV$ follows a Gaussian distribution with mean $\mean$
and covariance matrix $\cov_s$ parameterized by a positive
\emph{roughness parameter} $s$ that characterizes the roughness of the boundary.
While we focus on object boundaries, our work
applies to the general problem of estimating the roughness of a scalar function,
with many other applications as illustrated in~\Cref{sec:datafitting}.

With this form of the prior for $\bV$, we can think of $s$ as a
hyperparameter that adds uncertainty into the fractional differentiability of
the Whittle-Mat\'ern prior.
Thus, we can formulate a \emph{hierarchical posterior} of the form
  \begin{equation} \label{eq:posterior}
    \pi_{\bV\!, S | \bY }(\bv,s) \propto \pi_{\bY | \bV}(\by) \, \pi_{\bV | S}(\bv)
    \, \pi_S(s),
  \end{equation}
where $\pi_{\bY | \bV}(\by)$ is the \emph{likelihood} which depends on the forward model \eqref{eq:ip_deteministic} and the noise assumption.
The density $\pi_{\bV | S}(\bv) \pi_S(s)$ is the \emph{prior} in a hierarchical structure,
where $\pi_{\bV | S}(\bv)$ formulates the prior information on $\bv$ given the
roughness hyperparameter $s$, and $\pi_S(s)$ is the hyperprior for~$s$.
By means of this hierarchical posterior, we are able to infer the boundary and
the level of its roughness simultaneously, and we can quantify their uncertainties.

The resulting hierarchical posterior in \eqref{eq:posterior}
poses computational challenges due to the large-scale nature of the problem
and the lack of a closed-form expression for the solution.
In particular, the dependency of $\bV$ on $S$ poses challenges in exploring
the posterior using sampling methods \cite{mcbook}.
To alleviate the complexity of the sampling process, inspired by
the Karhunen-Lo\'eve (KL) expansion \cite{ibragimov2012gaussian},
we replace $\bV$ with a new random variable $\bU$ (see details in \Cref{sec:model})
that is independent of $S$.
Thus, we transform \eqref{eq:posterior} into the form
  \begin{equation}
    \pi_{\bU\!, S | \bY }(\bu,s) \propto \pi_{\bY | (\bU,S)}(\by) \, \pi_{\bU}(\bu)
    \, \pi_S(s) .
  \end{equation}
The independence of $\bU$ and $S$ allows us to develop efficient sampling methods.
Furthermore, by exploiting the relation between $\bU$ and $\bV$ we propose
efficient algorithms to compute the mapping between their realizations
based on the fast Fourier transform (FFT) algorithm.

%% file: 2.method.tex
\section{The Hierarchical Prior} \label{sec:model}

In this section, we elaborate on the hierarchical prior for $\bV$, which
depends on the roughness parameter $S$, and we introduce a variable transformation
to $\bU$ that does not depend on $S$.

Recall that $\bV$ follows a Whittle-Mat\'ern prior
in the form of $\bV \sim \mathcal{N} (\mean, \cov_s)$, where the
covariance matrix is parameterized by the roughness parameter~$s$.
In order to overcome the dependency of $\bV$ on $s$,
we introduce a mapping $\mathcal{F}$ to transform a standard normal multivariate random
variable $\bU \sim \mathcal N(\boldsymbol{0},\matI)$
to $\bV$.
This mapping overcomes the dependency between $\bV$ and $S$, and can be efficiently computed by means of the FFT algorithm, see \Cref{sec:method}.

We utilize the fact that the realization of the random function $\mathcal{V}$,
which characterizes the boundary, can be expanded in the trigonometric basis
  \begin{equation} \label{eq:basis}
    \left\{ \frac{1}{\sqrt{\pi}} \sin\left(j\frac{2\pi x}{\ell}\right),
    \frac{1}{\sqrt{\pi}} \cos \left( j\frac{2\pi x}{\ell} \right) \right\}_{j=0}^\infty,
  \end{equation}
where $x \in [0,\ell)$ and $\ell>0$ is the length of the periodic domain.
To implement this,
we introduce the \emph{truncation parameter} $k$ and only consider the basis functions
with indices $j=0,\dots,k$.
Note that these indices correspond to $2k+1$ basis functions since the
sine function for $j=0$ is the zero function.
Furthermore, we assume a \emph{deterministic} mean for
the random function $\mathcal{V}$, represented by a constant function.
Hence, it is natural to also omit the cosine basis function for $j=0$,
and thus we
expand $\mathcal V$ in the remaining $2k$ basis functions for $j=1,\ldots,k$.
The value of $k$ controls the number of frequencies that can be resolved
in a function expanded in the basis \eqref{eq:basis}, and therefore
we must choose $k$ large enough to allow the desired variations in $\mathcal{V}$.

We now discretize the $\sin$-$\cos$ basis functions in \eqref{eq:basis} on a uniform mesh
of size $m$, to obtain the vectors $\be_1,\be_n,\ldots,\be_{2k} \in \mathbb{R}^m$.
This provides a basis for the realizations of $\bV$.
Then we can factor $\cov_s$ as
  \begin{equation} \label{eq:cov}
    \cov_s = \matB \, \matL_s \, \matB^T,
  \end{equation}
where
  \begin{equation}
  \label{eq:B}
    \matB = [ \be_1, \be_2 , \ldots , \be_{2k} ]
    \in\mathbb{R}^{m \times (2k)}
  \end{equation}
and $\matL_s$ is a $(2k)$-by-$(2k)$ diagonal matrix
with pairwise identical diagonal elements $\lambda_1,\dots, \lambda_{2k}$ defined as
  \begin{equation}\label{eq:lambda}
    \lambda_{2j-1} = \lambda_{2j} = c_s(\sigma + j^2)^{-2(s + \nicefrac{1}{2})},
    \quad j=1, \dots, k.
  \end{equation}
Here, $\sigma>0$ is referred to as the \emph{length scale} and it controls the
correlation between the components of $\bV$.
Moreover, the constant $c_s$ is defined via
  \[
    c_s^{-1} = 2\sum_{j=1}^{\infty} (\sigma +j^2)^{-2(s + \nicefrac{1}{2})}
  \]
to ensure that $\sum_{j=1}^{\infty} \lambda_j = 1$;
it is often referred to as the global variance of Whittle-Mat\'ern covariance
\cite{1930-8337_2014_2_561}.
Note that $\lambda_j >0$ for all $j$ and that
$\cov_s$ is a symmetric semi-positive definite matrix
of rank $2k$
whose nonzero eigenvalues are proportional to $\lambda_j$.

According to the results in \cite{pham2021multivariate}, a Gaussian multivariate
random vector with the covariance matrix \eqref{eq:cov} can be expressed as
a linear combination of the $\sin$-$\cos$ basis:
  \begin{equation} \label{eq:kl_expansion}
    \bV =  \sum_{i=1}^{2k} U_i \sqrt{\lambda_j} \, \be_i,
  \end{equation}
where $\{U_i\}_{i=1}^{2k}$ are independent standard normal random variables.
Note that the variance for $\bV$ satisfies
  \[
    \mathrm{Var}(\bV) = \sum_{i=1}^{2k} \lambda_i \mathrm{Var}(U_i) \leq 1
  \]
and it tends to 1 as $k\to \infty$.

The expansion in \eqref{eq:kl_expansion} is called the \emph{Karhunen-Lo\'eve expansion}
in the limit where $k\to \infty$.
We refer the reader to \cite{ibragimov2012gaussian} for the notion of convergence
of the series as $k\to \infty$.
By defining the random vector $\bU = [U_1, \ldots, U_{2k}]^T$ and the mapping
$\mathcal{F}: \bU \mapsto \bV$, we can rewrite the inverse problem
\eqref{eq:ip_bayes} as
  \begin{equation} \label{eq:ip_goal}
    \bY = \mathcal G\circ \mathcal{F}(\bU) + \bE,
  \end{equation}
with the posterior
  \begin{equation} \label{eq:posterior_U}
    \pi_{\bU\!, S | \bY }(\bu,s) \propto \pi_{\bY|(\bU\!,S)}(\by) \,
    \pi_{\bU}(\bu) \, \pi_S(s).
  \end{equation}
Here, the likelihood density $\pi_{\mathbf Y|(\mathbf U,S)}$ follows the
noise model and the problem formulation \eqref{eq:ip_goal}.
Note that $\bU$ and $S$ are independent random variables
(while $\bV$ and $S$ are dependent in \eqref{eq:posterior}).
Therefore, the posterior \eqref{eq:posterior_U} allows a simpler sampling strategy
than the one for \eqref{eq:posterior}.

\begin{figure}
\centerline{\includegraphics[width=\columnwidth]{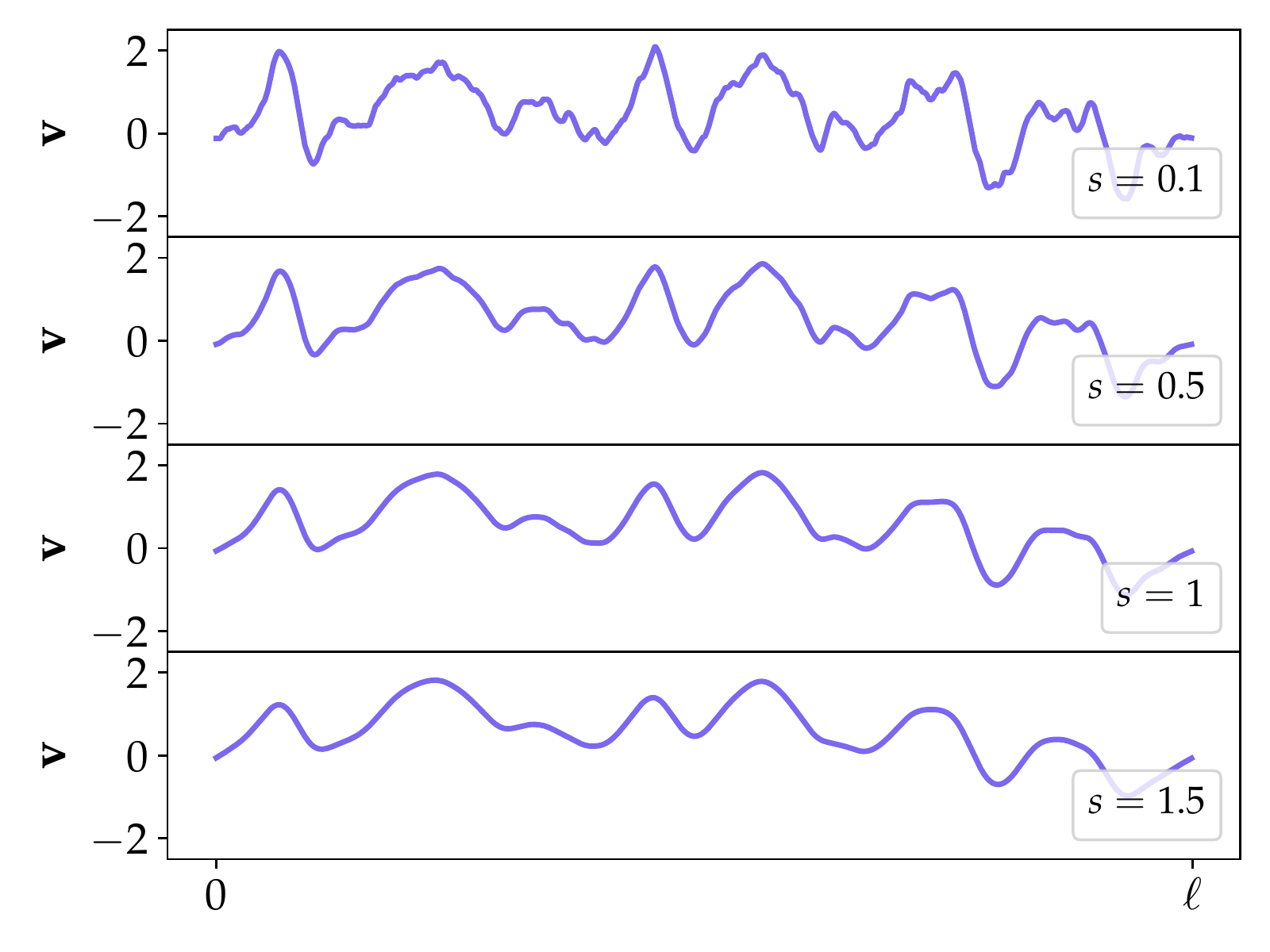} }
\caption{Samples from the random variable $\bV|S$ for a fixed realization
of $\mathbf u$ and four values of $s$.} \label{fig:prior}
\end{figure}

In the limit of $k\to \infty$, the parameter $s$ controls the level of differentiability/roughness of realizations of $\bV$.
In this case, realizations of $\bV$ almost surely belong to the fractional Hilbert space
of periodic functions $H_{\text{period}}^t$
for all $t < 2s + \nicefrac{1}{2}$ (cf.\ \cite{Dunlop2016}),
where the superscript $t$ indicates the fractional differentiability of
a periodic function (see \cite{DINEZZA2012521}).
The qualitative behavior of this smoothness is retained in the discretized case when the discretization mesh is fine enough, i.e., when $m$ is large enough.
\Cref{fig:prior} illustrates realizations of the prior distribution
$\pi_{\bV|S}(\bv)$ according to \eqref{eq:kl_expansion} with $m= 1024$,
for a fixed $\bu$ and different values of $s$.
We see that as $s$ increases the realizations of $\bV|S$ becomes ``smoother''
(less ``rough'').

\begin{remark}
In this work, $\bv$ is a discretization of a real-valued function.
Our findings can be naturally generalized to the case where $\bv$ is a discretization
of a function of $d$ variables
($\mathcal{V}:\mathbb R^d \to \mathbb R$ with $d\geq1$)
that is periodic in all $d$ dimensions.
In this case, we modify the values $\lambda_j$ to take the form
  \begin{equation} \label{eq:eigs_modified}
    \lambda_{2j-1} = \lambda_{2j} = c_s(\sigma + j^2)^{-2(s + \nicefrac{d}{2})},
    \quad j=1, \dots, k-1,
  \end{equation}
and the global variance is
  \begin{equation}
    c_s^{-1}= 2\sum_{j=1}^{\infty} (\sigma +j^2)^{-2(s + \nicefrac{d}{2})}.
  \end{equation}
This ensures that $\bV$ almost surely belong to the fractional Hilbert space
$H_{\text{$d$-{\rm period}}}^t$ for all $t < 2s + d/2$,
as $k\to \infty$, cf.\ \cite{Dunlop2016,DINEZZA2012521}.
\end{remark}

\section{Algorithms for the mappings $\mathcal{F}$ and $\mathcal{F}^{-1}$}\label{sec:method}

In this section we discuss how to compute the expansion \eqref{eq:kl_expansion}
for a realization $\bu$ without the need for storing the $\sin$-$\cos$ basis.
The key idea is to use the relation between this basis and the discrete
Fourier transform (DFT),
which allows us to use the FFT and IFFT algorithms for fast computations.

To describe the details, let $\boldf \in \mathbb{C}^{m}$ be a complex-valued vector given by
  \begin{equation} \label{eq:fourier}
    \boldf = \sum_{j=1}^{k} \hat{f}_j (\be_{2j} + i \be_{2j-1}  ),
  \end{equation}
where $i = \sqrt{-1}$ and $\{\hat f_j\in \mathbb C\}_{j=1}^{k}$
are expansion coefficients.
The real vectors $\be_{2j-1}$ and $\be_{2j}$ are
the columns of the matrix $\matB$ in \eqref{eq:B}, and hence
the complex vectors $\be_{2j} + i \be_{2j-1}$ are scaled
discretizations of the Fourier basis functions.
We note that $\boldf$ given by \eqref{eq:fourier} has zero mean.
The expansion coefficients $\hat{f}_1, \ldots , \hat{f}_{k}$ are the
coefficients of the DFT of $\boldf$, and we do not use
the coefficient $\hat{f}_0$ associated with the constant function.

At this time we choose to write the expansion coefficients in the form
  \begin{equation} \label{eq:fourier_expansion_coeffs}
    \hat{f}_j = \sqrt{\lambda_{2j}} \, ( g_j + i h_j ), \qquad g_j,h_j \in \mathbb R.
  \end{equation}
Using that $\lambda_{2j} = \lambda_{2j-1}$,
we can then rewrite \eqref{eq:fourier} to separate the real and the imaginary parts
of $\boldf$:
  \begin{equation} \label{eq:f}
    \begin{aligned}
        \boldf = \sum_{j=1}^{k} & \Big( \sqrt{\lambda_{2j}} \, ( g_j \mathbf e_{2j} -
          h_j \be_{2j-1} ) \\
        & + i  \sqrt{\lambda_{2j}} \, ( h_j \mathbf e_{2j} + g_j\mathbf e_{2j-1} ) \Big) .
    \end{aligned}
  \end{equation}
Then, by setting
  \begin{equation} \label{eq:gh}
    g_j = u_{2j}, \ h_j = u_{2j-1} , \qquad j=1,\ldots,k
  \end{equation}
we have that $\text{imag}(\boldf)$,
the imaginary part of the complex vector $\boldf$,
is in the form of a realization of \eqref{eq:kl_expansion}.

Now let the vectors $\bg , \bh \in \mathbb{R}^{k}$ contain $g_j$ and $h_j$,
respectively.
To ensure that the mapping $(\bg, \bh) \mapsto \text{imag}(\boldf)$
is invertible, we also need the real part of $\boldf$,
denoted by $\text{real}(\boldf)$.
We use the following relation \cite[\S 4]{doi:10.1137/1.9781611971514}:
  \begin{equation} \label{eq:v_sum}
  \begin{aligned}
    \frac{ \text{real}(\boldf) + \text{imag}(\boldf)}{\sqrt 2} =
    \sum_{j=1}^{k}& \Big( \sqrt{\lambda_{2j}} \, \frac{g_j + h_j}{\sqrt 2} \, \be_{2j} \\
    &+ \sqrt{\lambda_{2j}} \, \frac{g_j - h_j}{\sqrt 2} \, \be_{2j-1} \Big).
  \end{aligned}
  \end{equation}
According to the following proposition, the right-hand-side of \eqref{eq:v_sum} is
also in the form of a realization of \eqref{eq:kl_expansion}.
For a proof, see, e.g., \cite[\S 5.4]{Pitman}.

\begin{proposition}
Let $G$ and $H$ be independent and standard normal distributed real-valued random variables, then the random variables $(G+H)/\sqrt{2} $ and $(H - G)/\sqrt{2}$ are independent and standard normal distributed.
\end{proposition}

We can now construct the vector $\bv$ in terms of $\boldf$ as
  \begin{equation}
    \bv =  \frac{ \text{real}(\boldf) + \text{imag}(\boldf)}{\sqrt 2}.
  \end{equation}
Since the elements of $\boldf$ are related to those of $\bu$ via
\eqref{eq:f} and \eqref{eq:gh}, we have thus established a simple
computational relation between $\bu$ and $\bv$.

\Cref{alg:forward} summarizes the fast method to apply the mapping $\mathcal{F}$
to a random vector $\bu$, with the vectors
$\bg$ and $\bh$ being realizations from multivariate standard normal distributions.
Here, IFFT($\cdot$) denotes the fast inverse Fourier transform.

\begin{remark}
In many applications of FFT and IFFT it is common to let $m$ be a power of 2.
Since we use $2k$ expansion coefficients (because the coefficients for the
constant and zero functions are omitted), it is natural to use $m=2k$,
where $k$ is a power of 2.
FFT and IFFT require the same dimension $m$ for the input and output,
while we only have $k$ coefficients available in \eqref{eq:fourier_expansion_coeffs}.
We therefore need to supply additional coefficients (assumed to be zero)
associated with the frequencies in the set
  \begin{equation} \label{eq:Upsilon}
    \Upsilon = \{ j:j=0 \ \mathrm{or} \ k+1 \leq j \leq 2k \}
  \end{equation}
to obtain the vector
  \begin{equation} \label{eq:extend}
    \hat {\boldf} = [0,\hat{f}_1, \ldots, \hat{f}_{k}, 0, \ldots, 0]\in \mathbb R^{2k}.
  \end{equation}
Similarly, we can simply discard the components of $\hat{\boldf}$ corresponding to indices in
$\Upsilon$ to obtain the expansion coefficients $\hat{f}_1, \ldots, \hat{f}_{k}$.
\end{remark}

\begin{algorithm}
\caption{Compute $\bv = \mathcal{F}(\bu)$} \label{alg:forward}
\begin{algorithmic}[1]
\renewcommand{\algorithmicrequire}{\textbf{Input:}}
\renewcommand{\algorithmicensure}{\textbf{Output:}}
\Require $\bu \in \mathbb R^{2k}$ and coefficients $\{\lambda_j\}_{j=1}^{k}$.
    \State Construct $\bg, \bh \in \mathbb R^{k}$ from $\bu$ by \eqref{eq:gh}.
    \State Construct the $k\times k$ diagonal matrix
      $$\matL = \text{diag}(\sqrt{\lambda_1},\ldots, \sqrt{\lambda_{k}}).$$
    \State Set $\hat{\boldf}_{\mathrm{short}} = \matL ( \bg + i \bh)$
      with $\hat{\boldf}_{\mathrm{short}} \in \mathbb{C}^{k}$.
    \State Compute $\hat{\boldf} \in \mathbb C^{2k}$ by augmenting
      $\hat{\boldf}_{\mathrm{short}}$ with zeros according to \eqref{eq:extend}.
    \State Compute $\boldf = \mathrm{IFFT}(\hat{\boldf})$.
    \State Compute $\bv =  \bigl( \text{real}(\boldf) +
      \text{imag}(\boldf) \bigr) /\sqrt 2$.
\Ensure The vector $\bv$ in the $\sin$-$\cos$ basis.
\end{algorithmic}
\end{algorithm}

The computational complexity of \Cref{alg:forward} for
computing a sample from the distribution
$\mathcal N(\boldsymbol{0},\cov_s)$ for a given $s$ comprises computation of
the complex vector $\hat{\boldf}_{\mathrm{short}}$ in $\mathcal O(k)$ operations
and computing an IFFT of $\hat{\boldf}$ in $\mathcal O(2k\log(2k))$ operations.

To formulate an algorithm for the inverse mapping $\mathcal{F}^{-1}: \bv \mapsto \bu$
we assume that a vector $\bv$ is given in the form of \eqref{eq:v_sum}.
Then, we can compute a Euclidean projection of $\bv$
onto the Fourier basis $\{ \be_{2j} - i \be_{2j-1}\}_{j=1}^{k}$ by using that
  \begin{equation}
    \big\langle \bv, \be_{2j} - i \be_{2j-1} \big\rangle =
    \sqrt{\frac{\lambda_{2j}}{2}}\bigl( (g_j + h_j) + i (g_j - h_j) \bigr).
  \end{equation}
This is, in fact, the $j$th component of the DFT of $\bv$.
Therefore we can recover $g_j$ and $h_j$ from
  \begin{equation}
  \begin{aligned}
    g_j &= \sqrt{ \frac{2}{\lambda_{2j}} } \bigl( \text{real}(v_j) +
     \text{imag}(v_j) \bigr) \\
    h_j &= \sqrt{ \frac{2}{\lambda_{2j}} } \bigl( \text{real}(v_j) -
     \text{imag}(v_j) \bigr)
  \end{aligned}
  \end{equation}
for $j=1,\ldots, k$.
We summarize the method to apply the mapping $\mathcal F^{-1}$ in \Cref{alg:backward},
where FFT($\cdot$) denotes the fast Fourier transform.
The computational complexity of \Cref{alg:backward} is similar to that
of \Cref{alg:forward}.

\begin{algorithm}
\caption{Compute $\bu = \mathcal{F}^{-1}(\bv)$} \label{alg:backward}
\begin{algorithmic}[1]
\renewcommand{\algorithmicrequire}{\textbf{Input:}}
\renewcommand{\algorithmicensure}{\textbf{Output:}}
\Require $\bv \in \mathbb R^{m}$ and coefficients $\{\lambda_j\}_{j=1}^{k}$.
    \State Construct the $k \times k$ diagonal matrix
      $$\matL = \text{diag}(\sqrt{\lambda_1},\cdots, \sqrt{\lambda_{k}}).$$
    \State Compute the vector
      $\hat{\bv} = \mathrm{FFT}(\bv) \in \mathbb{C}^{2k}$.
    \State Compute $\hat{\bv}_{\mathrm{short}} \in \mathbb{C}^{k}$ by omitting
      the elements of $\hat{\bv}$ with indices in the set $\Upsilon$ \eqref{eq:Upsilon}.
    \State Compute $\bt = \matL^{-1} \hat{\bv}_{\mathrm{short}} \in \mathbb{C}^{k}$.
    \State Compute $\bg =\sqrt 2 \, (\text{real}(\bt) + \text{imag}(\bt))$.
    \State Compute $\bh =\sqrt 2 \, (\text{real}(\bt) - \text{imag}(\bt))$.
    \State Construct $\bu$ from $\bg$ and $\bh$ according to \eqref{eq:gh}.
\Ensure The vector $\bu$.
\end{algorithmic}
\end{algorithm}